\documentclass[12pt,letterpaper]{article}
\usepackage{epsfig}
\def\beq#1{\begin{equation}\label{#1}}
\def\eeq{\end{equation}}
\def\beqs#1{\begin{equation}\label{#1}\begin{array}{lllllllllllllllllll}}
\def\eeqs{\end{array}\end{equation}}
\def\beqss#1{\begin{equation}\label{#1}\left\{\begin{array}{lllllllllllllllllll}}
\def\eeqss{\end{array}\right.\end{equation}}
\def\eq#1{(\ref{#1})}

\def\defeq{\stackrel{\rm def}{=}}

\def\col{\mathop{\rm col} \nolimits}

\def\span#1{{\rm span}\left\{#1\right\}}

\def\dfrac#1#2{\displaystyle\frac{#1}{#2}}

\def\dlim{\displaystyle\lim}

\def\D#1#2{\displaystyle\frac{\partial#1}{\partial#2}}

\def\matrix#1{\left[\begin{array}{ccccccccccccccccccccc}#1\end{array}\right]}
 
\def\={&=&}

\renewcommand{\=}{=}

\newcommand{\cQ}{\mathcal{Q}}
\newcommand{\tx}{\tilde{x}}
\newcommand{\ts}{\tilde{s}}
\newcommand{\bs}{\bar{s}}
\newcommand{\bx}{\bar{x}}

\newcommand{\bv}{\bar{v}}

\newcommand{\bw}{\bar{w}}

\newcommand{\dx}{x^{\mbox{\rm \tiny D}}}
\newcommand{\dy}{y^{\mbox{\rm \tiny D}}}
\newcommand{\ds}{s^{\mbox{\rm \tiny D}}}
\newcommand{\du}{u^{\mbox{\rm \tiny D}}}
\newcommand{\dv}{v^{\mbox{\rm \tiny D}}}
\newcommand{\des}{{\mbox{\rm \tiny D}}}
\newcommand{\eproof}{\hfill \rule{1.5ex}{1.5ex}}
\newtheorem{theorem}{Theorem}
\newtheorem{proposition}{Proposition}
\newtheorem{remark}{Remark}
\newtheorem{definition}{Definition}
\begin{document}
\begin{center}
{\Large \bf Trajectory tracking control for} \\ {\Large \bf
maneuverable nonholonomic systems.}\\[10mm]
\end{center}
\begin{center}
{\large Sergei V. Gusev}\\
\end{center}
{\normalsize\it Department of Mathematics and Mechanics,}
{\normalsize\it St Petersburg State University } \newline
{\normalsize\it 2 Bibliotechnaya sq., Peterhof,
St.Petersburg, 198904, Russia}\newline
E-mail address: gusev@ieee.org
\begin{center}
{\large Igor A. Makarov}\\
\end{center}
{\normalsize \it
Institute for Problems of Mechanical Engineering,}
{\normalsize \it Russian
Academy of Sciences} \newline
{\normalsize \it 61 Bolshoy av., St.Petersburg, 199178,
Russia}\newline
E-mail address: mak@ccs1.ipme.ru\\
\begin{center} {\large \bf Abstract} \end{center}
The paper considers a motion control problem for kinematic models
of nonholonomic wheeled systems. The class of maneuverable wheeled
systems is defined consisting of systems that can follow any
sufficiently smooth non-stop trajectory on the plane. A sufficient
condition for maneuverability is obtained. The design of control
law that stabilizes motion along the desired trajectory on the
plane is performed in two steps. On the first step the trajectory
on the configuration manifold of the system and the input function
are constructed that ensure the exact reproduction of the desired
trajectory on the plane. The second step is the stabilization of
the constructed trajectory on the configuration manifold of the
system. For this purpose a recursive procedure is used that is a
version of backstepping algorithm meant for non-stationary systems
nonlinearly depending on input. The procedure results in the
continuous memoryless feedback that stabilizes the trajectory on
the configuration manifold of the system. As an example the motion
control problem for a truck with multiple trailers is considered.
It is shown that the proposed control law stabilizes the desired
trajectory of the vehicle on the plane for all initial states of
the system from some open dense submanifold of the configuration
manifold, i.~e., almost globally. The statement takes place both
for a truck pulling any number of trailers in a forward direction
and for a truck pushing any number of trailers in a backward
direction. The latter result is the solution of the intuitively
hard problem of the road train reverse motion control. The
effectiveness of the proposed control is demonstrated by
simulation. Animated examples are presented at Sergei~V.~Gusev Web
Page.\\[1em]



\newpage

\noindent
{\large \bf 1.~~Introduction.}\\

The control of nonholonomic mechanical systems
is a subject of intensive study
(see survey \cite{KMcC95}) that is mainly
devoted to the transport robot control.
These investigations can be classified
into two groups: feedforward  and feedback
control strategies.
The first direction, known as
the motion planning problem, is presented
in the comprehensive monograph \cite{Li&Ca92}.
The second direction can be subdivided into
problems of the equilibrium manifold stabilization
\cite{BRM92,MS94},
the zero state stabilization
\cite{BRM92,CdW&Sor92,C92,MCl&Mur97,Pom92,S95,Sor&Eg95},
and the stabilization of desired trajectory.
Despite that latter problem is of great
practical importance it is less investigated
then other stabilization problems.

One  approach to stabilization of the desired
trajectory uses  approximate linearization
of the system in the neighborhood of this
trajectory \cite{ABKSS93, WTSML94}.
Unfortunately, thus obtained linear feedback
is guaranteed to 
perform well only in a
small neighborhood of the desired trajectory.

Another approach is based on the system
transformation to the chained form \cite{S95}.
In \cite{JN99} such a transformation is used
for the trajectory stabilization of wheeled systems.
Though the constructed feedback is not globally stabilizing
its domain of attraction is not necessary a small
neighborhood of the desired trajectory as in
the previous case.
However, this approach has a drawback that
some trajectories, despite of being
traceable in principle, are
nevertheless missed in the
domain of the transformation and hence they
could not be stabilized.
Thus the approach imposes unnecessary and unnatural
restrictions on the desired trajectories
of the system.
For example, whether or not a particular trajectory
is stabilizable might depend on the choice
of a Cartesian coordinate system on the plane.

This paper deals with a specific variant of
a trajectory tracking problem for wheeled
systems.
We select a point fixed in the body of
the vehicle
and are trying to define the vehicle control
that stabilizes the motion of this
{\sl distinguished point}
along a given trajectory on the plane.
Note that when given a planar desired
trajectory,
we are not presented with
the corresponding high-dimensional
trajectory on the configuration
manifold of the system.
Therefore, we find the control in two steps:
1) the motion planning step during which
we construct an input function and
the corresponding trajectory
on the configuration manifold of the system,
so that the latter
maps exactly onto
the planar desired trajectory;
2) the stabilization step, wherein we
stabilize the
constructed trajectory.

The purpose of this paper is the design of the control law that
stabilizes any non-stop motion of the distinguished point of
the vehicle along any sufficiently smooth curve on the plane.
To this end it is supposed that the system can trace
any such a trajectory on the plane.
Systems that  have this property are called maneuverable.

We give a sufficient condition for
maneuverability of
wheeled systems. For systems that satisfy this condition the
motion planning problem is solved, i.e., the algorithm is proposed
that using the desired trajectory of the distinguished point constructs
the corresponding  trajectory on the configuration manifold of the
system.

The stabilization of the constructed trajectory
is based on a nonlinear
state feedback transformation
of the kinematic model of the system
to  a simplified cascaded form.
The recursive application of a backstepping-like procedure
to the cascaded system gives
a continuous memoryless  feedback
that stabilizes the trajectory.

In our previous papers it was shown that
a polar state transformation can be used
to achieve the local stabilization of the desired
trajectory in the case of caterpillar
\cite{G&M89} and four-wheeled \cite{M93b}
mobile robots.
Such a transformation allows to obtain the
high performance practical control algorithm
for trajectory tracking of the mobile platform
\cite{EHS01}.

This paper extends
previous results using a more general
state feedback transformation of the system.
The constructed  control law
stabilizes the motion of a
truck with multiple trailers  along any
non-stop trajectory
that has a sufficient number of bounded derivatives.
We tackle the task
both for a truck pulling any number of trailers
in a forward direction and even for a truck pushing
any number of trailers in a backward direction.
In addition, the stabilization is almost global
in sense that the  attraction domain of the trajectory
is an open dense submanifold of the
system configuration manifold.

For the Chaplygin sled\footnote{
The Chaplygin sled is the wheeled system that
has one wheel and two sleeve bearings \cite{NF72}.
The free movement of this system was first investigated
by Chaplygin in~\cite{Chaplygin}.
This system is also called knife edge \cite{BRM92}.},
which is the simplest wheeled system,
the constructed feedback is memoryless, continuous on the
whole configuration manifold and
guarantees the global stabilization of the desired
trajectory.
The latter example shows that
the assumption about the non-stop character of the motion
is essential. Because,  a well known Brockett's
result \cite{Br83} implies that a continuous memoryless feedback
cannot globally stabilize
the desired configuration of the Chaplygin sled.

The simulation shows the effectiveness of the proposed
method of trajectory tracking control.
The animated simulation is presented in~\cite{GWP}.
Preliminary results on the maneuverable vehicles
control can be found in~\cite{GM02}.

Finally, we should note that while the
present paper deals with  kinematic models
of wheeled systems,
the obtained results can serve as the basis for the
stabilization of dynamical models of vehicles
(by analogy
with results in \cite{GMPYL98},
where the stabilization of the dynamical model of a car
is considered)
as well as for the adaptive control of robots using
methods in
\cite{GMF91,GM94,GMPY99,GMPYL98}.
\\[1ex]

\noindent
{\large \bf 2.~~Mathematical model of the system.}\\[1ex]

Wheeled vehicles are nonholonomic mechanical systems.
We begin by describing the mathematical model of such systems.
The  configuration manifold of the system $\mathcal{Q}$
is the real smooth manifold with local coordinates
$q=(q_1,\ldots,q_N).$
Let $T_q \cQ$ and $T^\ast_q \cQ$ denote tangent and cotangent
spaces of the manifold $\mathcal{Q}$ at a point
$q,$ and let
$T \cQ = \bigcup_{q \in \cQ}T_q \cQ$ and
$T^\ast \cQ = \bigcup_{q \in \cQ}T^\ast_q \cQ$
denote tangent and cotangent bundles of this
manifold\footnote{Hereinafter we use conventional but
somewhat inaccurate notations that identify the point
of manifold and its coordinates.
Comments on
such shorthands
can be found in \cite{NS90}.}.
The kinematics of a nonholonomic system is described
by a set of one-forms
$\omega_i \in T^\ast \cQ, \; i = 1, \ldots,n,$
that define the
linear homogeneous nonholonomic constraints
\beq{NHC}
   \langle\omega_j(q), \dot{q}\rangle = 0,
   \ \ \ j=1,\ldots,n,
\eeq
where $ \langle\omega(q), \dot{q}\rangle$
denote the action of the linear functional
$\omega(q) \in T^\ast_q \cQ$ on the tangent vector
$\dot{q} \in T_q \cQ.$
The trajectory of a nonholonomic system is the function
$q \in C^1([0, \infty), \cQ)$
that satisfies the equations \eq{NHC}.
Hereinafter
$C^k(\mathcal{M}_1, \mathcal{M}_2), \; k=0, 1, 2, \ldots,$
denotes the class of $k$ times continuously differentiable
maps of the manifold $\mathcal{M}_1$
into the manifold $\mathcal{M}_2.$

Let $\mathcal{K}$ be  an open submanifold of $\cQ$
such that
the codistribution
$\Omega = \span{\omega_1, \ldots, \omega_n}
\subset T^\ast \cQ$ is constant-dimensional
on $\mathcal{K}$.
The codistribution $\Omega$ defines on $\mathcal{K}$
the $(m=N-n)$-dimensional  distribution
$\Delta = \Omega^{\bot} =
\{ g \in T\mathcal{K} \; | \;
\langle\omega, g\rangle \equiv 0 \;
\forall \omega \in \Omega \}.$
We assume that smooth vectorfields $g_1, \ldots, g_m$ form
a basis of the distribution $\Delta.$
Then any
lying in $\mathcal{K}$
trajectory of the nonholonomic system
satisfies the differential equation
\beq{KM}
   \dot{q} = \sum_{i=1}^m u_i g_i(q),
\end{equation}
where $u = \col(u_1, \ldots, u_m) \in
C([0, \infty), R^m).$
Equation \eq{KM} is referred to as a {\sl kinematic model}
of the nonholonomic mechanical system.
The freedom in defining the kinematic model
is in
the choice of the submanifold
$\mathcal{K}$ and the basis vectorfields
$g_1, \ldots, g_m.$
We consider \eq{KM} as the equation that describes
a control system with input $u.$


The {\sl wheeled system} is a set
of interconnected
rigid bodies with wheels that can move on the plane.
The wheels are constrained to roll without slipping.
The position of the system on the plane is defined by the
Cartesian coordinates $x_1, x_2$ of a distinguished point,
which is selected
belonging to one of the wheel axles.
Let $y_1$ be heading angle of the corresponding wheel,
then the
rolling without slipping constraint
for this wheel
takes the form
\beq{RNSC}
   \dot{x}_1\sin y_1-\dot{x}_2\cos y_1 = 0.
\eeq

We take a natural assumption that the constraint
equations are invariant with respect to translations
of the $x_1, x_2$ plane.

It is typical for a vehicular system
to have two scalar inputs.
In terms of nonholonomic constraints
it means that the number of constraints is two less
then the number of degrees of freedom of the system.

Our assumptions can be summed up
as follows:
\begin{enumerate}
\renewcommand{\theenumi}{\Roman{enumi}}
\item
$\mathcal{Q} = R^2 \times \check{\mathcal{Q}},$
where $\check{\mathcal{Q}}$ is a smooth manifold of
dimension $n.$ The vector of coordinates of the system
can be represented as $q = \col(x, y),$
where $x = \col(x_1, x_2)$ is the vector of
Cartesian coordinates of the distinguished point,
and $y = \col( y_1, \ldots y_n)$
are the remaining coordinates.
\item
The kinematics of the system is described by the
set of nonholonomic constraints \eq{NHC} that includes the
constraint \eq{RNSC}.
\item
The one-forms $\omega_i, i = 1, \ldots, n,$
do not depend on the coordinates $x_1, x_2.$
\end{enumerate}

Let $\mathcal{K}$ be an open submanifold of $\mathcal{Q}.$
It turns out that under very non-restrictive assumptions
the wheeled system admits on $\mathcal{K}$
the kinematic model of the following special type
\begin{eqnarray}
\label{CSx}
& & \begin{array}{lll}
   \dot{x}_1&=&u_1\cos y_1, \\
   \dot{x}_2&=&u_1\sin y_1,
\end{array} \\
\label{CSy}
& & \begin{array}{lll}
   \dot{y}    &=& u_1 h_1(y)+ u_2 h_2(y),
\end{array}
\end{eqnarray}
where $u_1$ and $u_2$ are inputs,
$u_1$ is the longitudinal velocity of the distinguished
point motion along the vector
$(\dot{x}_1,\dot{x}_2) \in R^2,$
that has slope $y_1,$
$h_1, h_2 \in T\check{\mathcal{Q}}.$
The output of the system 
is the distinguished point position $x = \col(x_1, x_2).$
Therefore the output trajectory of the system \eq{CSx}, \eq{CSy}
will be  referred to as the distinguished point trajectory.

The control system \eq{CSx}, \eq{CSy} is the subject
of investigation in this paper.

\newpage
\begin{proposition}\hspace{-1ex}{\bf :}
\label{Prop2}
Suppose 
the manifold $\mathcal{K}$
has the following properties:
\begin{itemize}
\item[K1.]
$\mathcal{K} = R^2 \times \mathcal{Y},$
where $\mathcal{Y}$ is an open submanifold of
$\check{\mathcal{Q}}.$
\item[K2.]
The dimension of the codistribution
$\Omega = \span{\omega_1, \ldots, \omega_n}$
is constant on $\mathcal{K}.$
\item[K3.]
For any
$q_0 \in \mathcal{K}$
there exist the trajectory of the system
$q \in C^1([-1, 1] \rightarrow \mathcal{K}),$
passing through the point
$q_0$ ($q(0) = q_0$),
and such that $\dot{x}(0) \not= 0.$
\end{itemize}
Then the wheeled system admits the kinematic model
\eq{CSx}, \eq{CSy} on the manifold $\mathcal{K}.$
\end{proposition}

The proof of Proposition~1 is given in Appendix~A.

\begin{figure}[t]
\begin{center}
\includegraphics[width=8.5cm]{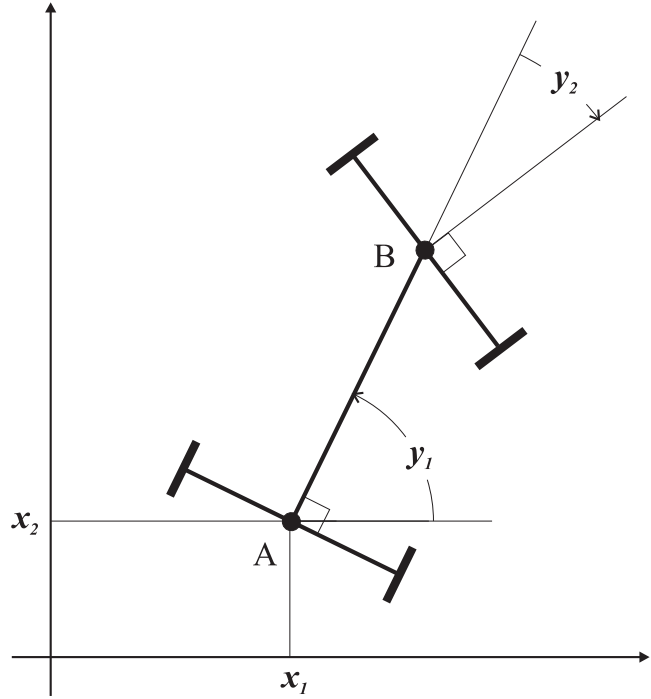}
\end{center}
\caption{Kinematic scheme of an automobile (I).}
\label{Fig1}
\end{figure}

As an example consider the kinematic model of an automobile, 
the scheme of that is shown in Fig.~\ref{Fig1}.
The vector of coordinates is $q = \col(x_1, x_2, y_1, y_2),$
where $x_1$ and $x_2$ are the Cartesian coordinates of
the midpoint of the rear axle,
$y_1$ is the heading angle, and $y_2$ is the
angle between the front and rear axles.
Define the configuration manifold as
$\mathcal{Q} = R^3\times S^1,$
where $S^1$ is the unit circle.
From here on we shall employ the  usual angular coordinate
on the $S^1$  and we shall define
trigonometric functions on $S^1,$ where the argument will be
the angle thus defined.

Nonholonomic constraints take the form \cite{MS93}
\beq{NCC}
\begin{array}{l}
\dot{x_1}\sin y_1 - \dot{x_2}\cos y_1 = 0,\\
\dot{x_1}\sin(y_1 + y_2) -
\dot{x_2}\cos(y_1 + y_2) - \dot{y_1}\cos y_2 = 0,
\end{array}
\end{equation}
where for simplicity the length of the  automobile base
(the segment AB)
is assumed to be unity.
The constraints \eq{NCC} are defined by the one-forms
$\omega_1(q) = \sin y_1 dx_1 - \cos y_1 dx_2, \;
\omega_2(q) = \sin(y_1 + y_2) dx_1 -
\cos(y_1+y_2) dx_2 - \cos y_2 dy_1.$
It is easy to see that
$\mathcal{K} = \{ q \in \mathcal{Q} \; | \;
\cos y_2 \not= 0 \}$ is the maximal submanifold
of the manifold $\mathcal{Q},$ where the codistribution
$\Omega = \span{\omega_1, \omega_2}$
is constant-dimensional.
The defined on $\mathcal{K}$ kinematic model
of an automobile has the known form \cite{MS93}
\beq{3wKM}
\begin{array}{l}
\dot x_1=u_1\cos y_1,\\
\dot x_2=u_1\sin y_1,\\
\dot y_1=u_1\tan y_2,\\
\dot y_2=u_2,
\end{array}
\eeq
where $u_1$ is the longitudinal velocity of the point A
and $u_2$ is the angular velocity of the front axle spin
relative to the automobile body.

Let us split the manifold $\mathcal{K}$ as
$\mathcal{K}=R^2\times \mathcal{Y},$
where
\beq{Y3wKM}
\mathcal{Y}=\{ y = (y_1, y_2) \in R\times S^1 \; |
\; \cos y_2 \not= 0 \}.
\eeq
Then the equation \eq{3wKM} takes the form \eq{CSx}, \eq{CSy}
with $h_1(y) = \col(\tan y_2, 0), \;
 h_2(y) = \col(0, 1).$
\\[1ex]

\noindent
{\large \bf 3.~~Maneuverable systems.}\\[1ex]

Let us consider the planar trajectory $\dx \in C([0,\infty), R^2)$
as the desired trajectory of the distinguished
point of the system \eq{CSx}, \eq{CSy}.
\begin{definition}\hspace{-1ex}{\bf :}
The trajectory $\dx$ is called admissible
trajectory of the distinguished
point of the system \eq{CSx}, \eq{CSy}
if $\dx \in C^{n+1}([0,\infty), R^2)$ and
\begin{equation}
\inf_{t \ge 0}|\dot{x}^\des (t)| > 0.
\label{|xd'|>0}
\end{equation}
\end{definition}
The set of all admissible trajectories of the distinguished
point of the system \eq{CSx}, \eq{CSy}
is denoted $X \subset C^{n+1}([0,\infty), R^2).$
\begin{definition}\hspace{-1ex}{\bf :}
The system \eq{CSx}, \eq{CSy}
is called maneuverable on an open submanifold
$\mathcal{M} \subset \mathcal{K},$
if for any admissible trajectory of the distinguished point
of the system there are the trajectory
$q^\des  = \col(x^\des, y^\des)
\in C^1([0, \infty), \mathcal{M})$
and the input
$\du \in C([0, \infty), R^2)$
that satisfy  \eq{CSx}, \eq{CSy}.
The operator
$M : X \to C^1([0, \infty), \mathcal{M})\times
C([0, \infty), R^2),$
which takes $\dx$ to the pair $(q^\des, u^\des),$
is called the maneuvering operator of the system.
When $\mathcal{M} = \mathcal{K},$
the system is called maneuverable (without specifying 
the manifold). 
\end{definition}

A wheeled system is not necessary maneuverable.
We shall demonstrate this with an example
at the end of the section.
The theorem below gives a sufficient condition
for the wheeled system maneuverability.
The proof of the theorem constructively defines the
set of maneuvering operators for the system.

Let us introduce some notation.
Consider a vector field $h$ and a function $\phi$  defined
on a manifold $\mathcal{Y}.$
Denote
$L_{h}\phi = \langle{\rm d}\phi, h\rangle =
\sum_{i=1}^k \frac{\partial \phi}{\partial y_i} h^{(i)}$
the Lie derivative of the function $\phi$ along  the
vector field $h.$
The repeated Lie derivatives
$ L^i_{h}\phi, i = 0, 1, 2, \ldots,$
are inductively defined by
$L^0_{h}\phi = \phi$ and $ L^i_{h}\phi =
L_{h}L^{i-1}_{h}\phi, \; i \ge 1.$
In what follows we suppose that all manifolds, vectorfields, and
functions are smooth enough to define all necessary Lie
derivatives.
\newpage
\begin{theorem}\hspace{-1ex}{\bf :}
\label{maneuver}
Let
$\mathcal{O} \subset \mathcal{Y}$
be an open submanifold 
that is split as $\mathcal{O} = R\times \check{\mathcal{O}}.$
Suppose that for all $y \in \mathcal{O}$
the following conditions hold
\begin{equation}
L_{h_2}L_{h_1}^i y_1 = 0, \;
i = 0, \ldots, n-2,
\label{Liy=0}
\end{equation}
\begin{equation}
L_{h_2}L_{h_1}^{n-1} y_1 \not= 0.
\label{Lky/=0}
\end{equation}
Then the smooth change of coordinates
\beq{S}
s = S(y),
\eeq
where the map $S \in C^1(\mathcal{O}, R^n)$
is defined by the equation
\begin{equation}\label{S=}
S(y)=\col(L_{h_1}^{0} y, \ldots, L_{h_1}^{n-1} y),
\end{equation}
and the nonsingular feedback transformation
\beq{F}
v = F(y)u,
\eeq
where the map $F \in C(\mathcal{O}, GL(2))$
is defined by the equation
\beqs{F=}
F(y) = \matrix{1, & 0\\
L_{h_1}^n\phi(y), & L_{h_2}L_{h_1}^{n-1}\phi(y)},
\eeqs
transforms the system \eq{CSx}, \eq{CSy} into the system
\begin{eqnarray} 
\label{xCS}
& & \begin{array}{l}
    \dot{x}_1=v_1 \cos s_1, \\
    \dot{x}_2=v_1 \sin s_1,
\end{array}\\
\label{sCS}
& & \begin{array}{l}
    \dot{s}_i=v_1 s_{i+1}, \; i = 1, \ldots, n-1, \\
    \dot{s}_n=v_2.
\end{array}
\end{eqnarray}
If $S$ bijectively maps $\mathcal{O}$ on $R^{n},$
then the system \eq{CSx}, \eq{CSy} is maneuverable on the
manifold
$\mathcal{M} = R^2\times \mathcal{O}.$
\end{theorem}

\begin{remark}\hspace{-1ex}{\bf :}
Under the transformations
\eq{S} and \eq{F} the equations
$s_1=y_1, v_1=u_1$ hold, and, consequently,
the $x$-subsystem \eq{CSx} does not change.
\end{remark}

{\bf Proof.}
The representation \eq{xCS}, \eq{sCS}
can be obtained using the well known
transformation of a nonlinear system to the
canonical linear one~\cite{Isidori}.
However it is not  difficult to prove this statement  directly.
Suppose that $y$ is the solution of the system \eq{CSy}
that corresponds to the input $u,$
and that $s$ and $v$ are defined by the transformations
\eq{S} and \eq{F} respectively.
Then from \eq{Liy=0}, \eq{S=},  and \eq{F=}
we obtain
\begin{eqnarray*}
& & \begin{array}{l}
\dot{s}_i(t) = u_1(t) L^i_{h_1}y_1(t) +
u_2(t) L_{h_2}L^{i-1}_{h_1} y_1(t) =
v_1(t) s_{i+1}(t), \;
i=1, \ldots, n-1,
\end{array}
\\
& & \begin{array}{l}
\dot{s}_n(t) = u_1(t) L^n_{h_1}y_1(t) +
u_2(t) L_{h_2}L^{n-1}_{h_1} y_1(t) = v_2(t).
\end{array}
\end{eqnarray*}
It is easy to see that, conversely, for any
solution of \eq{sCS} there corresponds some
solution of \eq{CSy}, which is defined  by the reverse
transformation of variables.

Using  \eq{xCS}, define the longitudinal
velocity of the distinguished point that moves
along the desired trajectory $\dx$
\beq{RefU1}
\dv_1(t) =
\pm \sqrt{(\dot{x}^\des_1(t))^2+(\dot{x}^\des_2(t))^2},
\eeq
where the sign of $\dv_1$ can be chosen arbitrary
but does not vary in time.
Calculating $\dot{s}_1(t)$ by virtue of
\eq{xCS} with $v_1(t) = \dv_1(t),$
we get
\beq{dy1}
    \dot{s}_1(t)=
    \dfrac{\ddot{x}^\des_2(t)\cos s_1 -
    \ddot{x}^\des_1(t)\sin s_1}{\dv_1(t)}.
\eeq
Note that \eq{|xd'|>0}
implies the inequality $\inf_{t \ge 0}|\dv_1(t)|>0.$
The desired heading angle
$\ds_1(t), \; t \ge 0,$
can be found as the solution of the differential equation
\eq{dy1}
with the initial value
$\ds_1(0)$
that satisfies the equation
\beq{s1(0)}
\dot{x}^\des(0) = \dv_1(0)
\matrix{\cos \ds_1(0)\\ \sin \ds_1(0)}.
\eeq
Define  $\ds(t)$ using the recursive formulae
\beq{ds_i=}
\ds_i(t)=\dot{s}^\des_{i-1}(t)/\dv_1(t),\
i=2,\ldots,n,
\eeq
and put
\beq{dv_2=}
\dv_2(t)=\dot{s}^\des_n(t).
\eeq
Then the triplet $\dx, \ds, \dv$ satisfies
the system of differential equations \eq{xCS}, \eq{sCS}.
Since the transformation $S$ is diffeomorphism
$\mathcal{O}$ onto $R^{n},$ we can define the trajectory
\beq{dy=}
\dy(t)=S^{-1}(\ds(t)),
\eeq
which satisfies the inclusion
$\dy(t)\in \mathcal{O}$ for all $t \ge 0.$
The desired input $\du$ can be uniquely
defined from the equation
\beq{du=}
\du(t) = F^{-1}(\dy(t)) \dv(t),
\eeq
because the matrix $F(y)$ is nonsingular in $\mathcal{O}.$
The triplet $\dx, \dy, \du$ satisfies \eq{CSx}, \eq{CSy}.
\eproof

The proof gives the procedure for determining
the maneuvering operators for the system
\eq{CSx}, \eq{CSy}.
After choosing the sign of $v^\des_1$ in 
\eq{RefU1} and the initial value $s^\des_1(0)$ satisfying
\eq{s1(0)}, the formulae
\eq{RefU1} -- \eq{du=}
uniquely define the pair $q^\des = \col(\dx, \dy), \du$
for a given admissible trajectory $\dx.$

As an example of Theorem~\ref{maneuver} application
let us show that the kinematic model of an automobile \eq{3wKM}
is maneuverable.
The manifold $\mathcal{Y}$ defined by \eq{Y3wKM}
is disconnected and consists of two connected components.
It is easy to see that the  conditions
\eq{Liy=0} and \eq{Lky/=0}
of Theorem~\ref{maneuver} are fulfilled
on each  component.
For definiteness choose one component
$\mathcal{O} = \{ (y_1, y_2) \in R^2
\; | \; |y_2| < \pi/2 \}.$
On $\mathcal{O}$ the formulae \eq{S=} and \eq{F=} define
the state transformation
$S(y) = \col(y_1, \tan y_2)$
and the feedback transformation
$v_1 = u_1, \; v_2 = u_2 \cos^{-2} y_2.$
The system \eq{3wKM} is maneuverable on the manifold
$\mathcal{O}\times R^2$ because $S$ bijectively maps
$\mathcal{O}$ onto $R^2.$

In this example the choice of sign "+" in \eq{RefU1}
defines the trajectory $q^\des$ and the input $\du$ that correspond
to an automobile forward motion along the desired
trajectory $\dx,$
whereas
the sign "-" in \eq{RefU1}
corresponds to an automobile
backward motion along the same trajectory $\dx.$

\begin{figure}[t]
\begin{center}
\includegraphics[width=8.5cm]{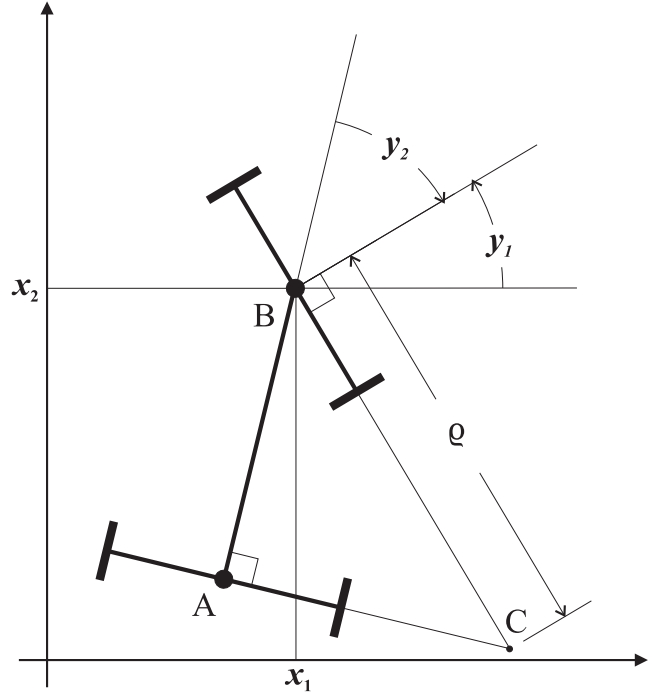}
\end{center}
\caption{Kinematic scheme of an automobile (II).}
\label{Fig2}
\end{figure}

Our next example 
shows that
not any wheeled system is maneuverable and,
more over, the maneuverability of the system depends on
the choice of coordinates.
Let us define the coordinates of automobile
as follows:
$x_1,x_2$ are the Cartesian coordinates of the distinguished point
B, which is the midpoint
of the front axle, $y_1$ is the heading angle of the front wheels,
$y_2$ is the angle between the front and rear axles
(see Fig.~\ref{Fig2}).
The configuration manifold of the system is
$\mathcal{Q} = R^3\times S^1$
and can be split as
$\mathcal{Q} = R^2 \times \mathcal{Y},$
where
$\mathcal{Y} = R \times S^1.$
The kinematic model is defined on the whole manifold
$\mathcal{Q}$ and takes the form \eq{CSx}, \eq{CSy}
with $h_1(y) = \col(\sin y_2, 0), \;
h_2(y) = \col(1, 1).$
It is easy to see that
conditions of Theorem~\ref{maneuver}
are not fulfilled for this system.

In fact, this system is not maneuverable.
To show this let us construct the
trajectory $\dx$ that cannot
be traced by the system.
From Fig.~\ref{Fig2} it is clear that
at each instant of time  $t$
the curvature of the point B trajectory
is equal to $1/\rho(t),$ where $\rho(t)$ is the length
of the hypotenuse of triangle ABC.
Hence, $\rho(t)\ge 1$ for all $t \ge 0.$
It follows 
that the system cannot trace any trajectory
$\dx$ the  curvature of which is less than one at some
instant of time $t \ge 0.$
Thus the considered system is not maneuverable.\\[1ex]

\noindent
{\large \bf 4.~~Formulation of the trajectory
stabilization problem.}\\[1ex]

Let $\dx$ be a desired trajectory of the distinguished
point of the system.
The  problem under consideration
is to stabilize the motion of
the distinguished point with
coordinates $x$
along the desired trajectory $\dx.$
But it is important also to guarantee boundedness
of the control input of the system.
Otherwise the stabilization has no
practical application.
For the input to be bounded,
it does not suffice
to assume the admissibility of the trajectory;
more strict requirements have to be
satisfied.
\newpage
\begin{definition}\hspace{-1ex}{\bf :}
The trajectory $\dx$ of the distinguished point
of the system \eq{CSx}, \eq{CSy}
is called strongly admissible if it is admissible
and, in addition, all its derivatives up to the order
$n+1$ are bounded, i.~e.,
$
\sup_{t \ge 0} \max_{k = 1, \ldots, n+1}|d^k/dt^k \dx(t)|
< +\infty.
$
\end{definition}

The set of strongly admissible trajectories
of the distinguished point
of the system \eq{CSx}, \eq{CSy} is denoted $\bar{X}.$

For maneuverable systems
the problem under consideration can be refined.
Let the system \eq{CSx}, \eq{CSy} be maneuverable
on the manifold $\mathcal{M} = R^2\times \mathcal{O}.$
Then for any admissible trajectory
on the plane
$\dx$
the maneuvering operator $M$ defines
the trajectory $q^\des = \col(\dx, \dy)$
of the system \eq{CSx}, \eq{CSy}
and the input $\du$ such
that the trajectory lies in 
$\mathcal{M}$
and corresponds to the input $\du.$
In this case becomes possible to replace
stabilization of
the desired planar trajectory $x^\des$
with that of the trajectory $q^\des.$
Clearly, the solution of the latter, more general problem, implies
the solution for the problem of tracking the trajectory $\dx.$

Suppose that $d$ is a metric on
the manifold $\mathcal{O}.$
Consider a map
$\Phi \in
C(\mathcal{M}\times \mathcal{M}\times \mathcal{U}, R^2),$
where
$\mathcal{U} = \{ u \in R^2 \; | \; u_1 \not= 0 \}.$
\begin{definition}\hspace{-1ex}{\bf :}
The feedback
\beq{u=Phi}
u = \Phi(q, q^\des(t), \du(t)),
\eeq
where $q = \col(x, y), \; q^\des = \col(\dx, \dy),$
stabilizes the trajectory $q^\des$ in $\mathcal{M},$
if for any initial value
$q(0) \in \mathcal{M}$
the solution of the closed-loop system
\eq{CSx}, \eq{CSy}, \eq{u=Phi}
is defined and lies in $\mathcal{M}$ for all $t \ge 0,$
and  if the following limits hold: 
\begin{eqnarray}
\label{x-dx->0}
\dlim_{t \to \infty}(x(t) - \dx(t)) = 0,\\
\label{y-dy->0}
\dlim_{t \to \infty}d(y(t), \dy(t)) = 0,\\
\label{u-du->0}
\dlim_{t \to \infty}(u(t) - \du(t)) = 0.
\end{eqnarray}
\end{definition}

This definition implies, in particular,
that if the input function
$\du$ is bounded, so will be the input $u$
of the closed-loop system.

Now we return to the problem of the desired trajectory $\dx$
stabilization.
Consider an operator
$U : [0,+\infty)\times \mathcal{M}\times \bar{X} \to
\mathcal{U},$ which for
$t \in [0,+\infty), \; q \in \mathcal{M},$ and
$\dx \in \bar{X}$ is defined as the superposition
\beq{U=1}
U(t, q, \dx) = \Phi(q, q^\des(t), \du(t)),
\eeq
where
\beq{U=2}
(q^\des, \du) = M(\dx).
\eeq
\begin{definition}\hspace{-1ex}{\bf :}
We say that the control law
\beq{u=U}
u = U(t, q, \dx)
\eeq
where the operator $U$ is defined by \eq{U=1}, \eq{U=2},
solves the problem of
stabilizing
the distinguished point trajectories
of the system \eq{CSx}, \eq{CSy}
on the manifold $\mathcal{M}$,
if for any strongly admissible trajectory $\dx$
the feedback \eq{u=Phi} stabilizes
the trajectory $q^\des$ defined by \eq{U=2}.
\end{definition}

The prime objective  of the paper is
solving the trajectories stabilization problem for
the wheeled  system \eq{CSx}, \eq{CSy} on the
manifold where the system is maneuverable.
To do this it is necessary to design
two operators $M$ and $\Phi.$
The former is already defined by
Theorem~\ref{maneuver}. The latter
is constructed in Section~5.

An additional objective is to make
as wide as possible
the domain where
the constructed control law solves
the trajectories stabilization problem.
From the practical standpoint it is desirable
to design the control law that solves this problem
on the whole configuration manifold of the system
$\mathcal{Q},$ i.e., globally.
Such a control law for the Chaplygin sled is described
below, in Section~6.
The general problem of global stabilization
of trajectories is not solved in the paper.
However, we believe our result closely
approximates the goal of the global stabilization.
Let us explain in what sense.

Suppose that
$\mathcal{K} = \bigcup_{i = 1}^m \mathcal{M}_i,$
where $\mathcal{M}_i, \; i = 1, \dots, m,$ are disjoint
open connected components of $\mathcal{K},$
and for every $i = 1, \dots, m,$  an operator
$U_i : [0,+\infty)\times \mathcal{M}_i\times \bar{X} \to
\mathcal{U}$ is defined.
Consider an operator
$U : [0,+\infty)\times \mathcal{K}\times \bar{X} \to
\mathcal{U}$ defined for
$t \in [0,+\infty), q \in \mathcal{K},
\dx \in \bar{X}$
by the equation
\beq{U=ifUi}
U(t, q, \dx) =  U_i(t, q, \dx),
\mbox{\rm ~if~ } q \in \mathcal{M}_i.
\eeq
\begin{definition}\hspace{-1ex}{\bf :}
We say that
the control law \eq{u=U} with the operator \eq{U=ifUi}
solves the problem of almost global stabilization
of the distinguished point trajectories of the system
\eq{CSx}, \eq{CSy}, if the manifold $\mathcal{K}$
is dense in the configuration manifold $\mathcal{Q},$
and if for every  $i = 1, \dots, m,$
the control law \eq{u=U} with $U=U_i,$
solves the problem of stabilizing
the distinguished point trajectories
of the system \eq{CSx}, \eq{CSy}
on the manifold $\mathcal{M}_i.$
\end{definition}

In the sense of this definition we shall show that
the proposed below control law solves  the problem of
almost global stabilization of the distinguished point
trajectories for the considered in
Section~6 kinematic model of
a truck with multiple trailers.\\

\noindent
{\large \bf 5.~~Stabilization of trajectories.}\\

Suppose that the system \eq{CSx}, \eq{CSy}
satisfies the conditions of Theorem~\ref{maneuver}
on the manifold $\mathcal{M}$
and that $M$ is a corresponding maneuvering
operator.
Let $\dx$ be an admissible trajectory
of the distinguished point.
Using the operator $M,$  define
the trajectory $q^\des = \col(\dx, \dy)$ and
the input $\du$ ($(q^\des, \du) = M(\dx)$).
After the state feedback transformation \eq{S}, \eq{F}
the system \eq{CSx}, \eq{CSy} takes on the
cascaded form \eq{xCS}, \eq{sCS} and our design
of the stabilizing feedback is based on that form.
Using the state transformation \eq{S} define
the trajectory
$s^\des(t) = S(y^\des(t)),$
$s^\des \in C^1([0,+\infty), R^n).$
To apply the backstepping technics
we shall represent the
system \eq{xCS}, \eq{sCS} in terms of 
deviations  from the trajectory $(\dx, s^\des).$\\

\noindent
{\bf 5.1.~Further transformation of the kinematic model.}
Consider the function
$
\tau(t)=\int_0^t\sqrt{(\dot{x}^\des_1(\kappa))^2
+ (\dot{x}^\des_2(\kappa))^2}d\kappa;
$
the value $\tau(t)$ is the length of the path traveled
by the distinguished point in a time $t.$
Due to the equality
$\dot\tau(t)=
|\dot{x}^\des(t)|$
and the inequality
\eq{|xd'|>0},
the function $\tau(.)$ maps
the interval $[0, \infty)$ bijectively onto itself.
The inverse function is denoted $t(.).$

Define functions
$\bx(\tau) = \dx(t(\tau)),
\bs(\tau) = \ds(t(\tau)),
\bv(\tau) = \dv(t(\tau)),
\tx(\tau) = x(t(\tau)) - \dx(t(\tau)),$
and
$\ts(\tau) = s(t(\tau)) - \ds(t(\tau)).$
The change of variables
from $t, x, s$ to $\tau, \tx, \ts$
transforms the system \eq{xCS}, \eq{sCS} into
\begin{eqnarray}
\label{x-CS}
\begin{array}{rcl}
    \tx' & = &
    - \bw_1 \matrix{\cos \bs_1 \\ \sin \bs_1} +
    w_1 G(\bs_1) \matrix{\cos \ts_1 \\ \sin \ts_1},
\end{array}\\
\label{s-CS}
\begin{array}{rcl}
    \ts'_i & = & w_1\ts_{i+1}+(w_1-\bw_1)\bs_{i+1},
    \indent i=1,\ldots,n-1, \\
    \ts'_n & = & w_2,
\end{array}
\end{eqnarray}
where
$$
   G(\bs_1) = \matrix{\cos \bs_1 & -\sin \bs_1 \\
   \sin \bs_1 & \cos \bs_1},
$$
\beq{w}
w_1(\tau) = v_1(t(\tau))/|\bv_1(\tau)|, \;
w_2(\tau) = (v_2(t(\tau)) - \bv_2(\tau))/|\bv_1(\tau)|.
\eeq
Note that 
$\bw_1 = \bv_1/|\bv_1| = \mbox{\rm sign } \bv_1$
does not depend on $\tau.$
Hereafter the prime denotes the differentiation
with respect to $\tau.$
The above transformation uses the following formulae:
$\bv_1(\tau) = \du_1(t(\tau)),$
$$
G(\bs_1) \matrix{\cos \ts_1 \\ \sin \ts_1} =
\matrix{\cos(\ts_1 + \bs_1) \\ \sin(\ts_1 + \bs_1)},
$$
\begin{eqnarray}
\label{x-CSdes}
& & \begin{array}{lcl}
    {\bx}'_1&=&\bw_1 \cos \bs_1, \\
    {x}'_2&=&\bw_1 \sin \bs_1,
\end{array}
\\
\label{s-CSdes}
& & \begin{array}{lcl}
    {\bs}'_i&=&\bw_1 \bs_{i+1}, \indent i=1,\ldots,n-1,\\
    {\bs}'_n&=&\bv_2/|\bv_1|.
\end{array}
\end{eqnarray}\\

\noindent
{\bf 5.2.~Cascaded system stabilization theorem.}
To design a stabilizing feedback for  the system
\eq{x-CS}, \eq{s-CS}
we use a recursive procedure
based on the idea of backstepping~\cite{KKK95}.
Let us formulate one step of the procedure.

Consider the cascaded system
\begin{eqnarray}
\label{BSCSz}
    z' = B(z, \zeta, p(\tau)),\\
\label{BSCSzeta}
    \zeta'=b(z, \zeta, p(\tau)) +
    \beta(z, \zeta, p(\tau)) \upsilon,
\end{eqnarray}
where
$z \in R^{k_z}, \zeta, \upsilon \in R, \tau \ge 0,$
$p \in C^1([0, \infty), R^{k_p}),$
$B,  \D{B}{\zeta} \in
C(R^{k_z}\times R\times R^{k_p}, R^{k_z}),$
$ b, \beta \in C(R^{k_z}\times R\times R^{k_p}, R),$
$\beta(z, \zeta, p(\tau) ) \ne 0$ for all
$z, \zeta, \tau.$

Suppose functions
$\alpha \in C^1(R^{k_z}\times R^{k_r}, R)$
and
$V \in C^1(R^{k_z}\times R^{k_p}\times R^{k_r}, R)$
are given that satisfy  the conditions
\beq{alpha=0}
\forall r
\;\; \alpha(0, r) = 0,
\eeq
\beq{V=0}
\forall \; p, r
\;\; V(0, p, r) = 0,
\eeq
\beq{V>0}
\forall \; z\not= 0, p, r
\;\; V(z, p, r) > 0,
\eeq
\beq{infV>0}
\forall \varepsilon>0, \mathcal{E}>0 \;
\inf_{|z| > \varepsilon, \; \max(|p|,|r|) < \mathcal{E}}
V(z, p, r) > 0.
\eeq
Define a function
$\hat{\alpha} :
R^{k_z+1} \times R^{k_p} \times R^{k_r}\times R^{k_r}
\to R$
by
\beq{alpha^}
 \begin{array}{ll}
    \hat{\alpha}(\hat{z}, p, r, r_1) = &
    {\beta(z, \zeta, p)^{-1}}
    \left[
    \D{\alpha(z, r)}{z} B(z, \zeta, p)
    + \D{\alpha(z, r)}{r} r_1 - \right.\\
    & \left.
    \delta^{-1} \D{V(z, p, r)}{z}
    D(z, \zeta, \alpha(z, r), p) -
    b(z, \zeta, p) -
    \gamma(\zeta - \alpha(z, r))
    \right],
 \end{array}
\eeq
where
$$
D(z, \zeta, \eta, p) =
\left\{
\begin{array}{ll}
\dfrac{B(z, \zeta, p)-B(z, \eta, p)}
{\zeta-\eta}, &
\mbox{\rm if } \zeta \not= \eta,\\
\D{B(z, \zeta, p)}{\zeta}, &
\mbox{\rm if } \zeta = \eta,
\end{array}
\right.
$$
and a function
$\hat{V} :
R^{k_z+1}\times R^{k_p} \times R^{k_r} \to R$
by
\beq{V^}
\hat{V}(\hat{z}, p, r) =
V(z, p, r)+ \delta (\zeta - \alpha(z, r))^2/2.
\eeq
Here $\gamma > 0, \; \delta > 0$ are parameters,
$z \in R^{k_z}, \; \zeta, \eta \in R, \;
\hat{z} = \col(z, \zeta), \; \;  p \in R^{k_p}, \;
r, r_1 \in R^{k_r}.$

\begin{theorem}\hspace{-1ex}{\bf :}
\label{backstep}
The function $\hat{\alpha}$ is continuous,
the function $\hat{V}$ is differentiable and
satisfies  the conditions \eq{V=0}--\eq{infV>0},
where $z$ should be replaced with $\hat{z}.$
\begin{itemize}
\item
If the  derivative of the function
$V(z(\tau), p(\tau), r(\tau))$
along the trajectories of the closed-loop system \eq{BSCSz},
\beq{zeta=alpha}
\zeta = \alpha(z, r(\tau))
\eeq
satisfies the inequality
\beq{V'<}
    V'(z(\tau), p(\tau), r(\tau))
    \le -2 \gamma V(z(\tau), p(\tau), r(\tau)),
\eeq
then the derivative of the function
$\hat{V}(z(\tau), p(\tau), r(\tau))$
along the trajectories of the closed-loop system
\eq{BSCSz}, \eq{BSCSzeta},
\beq{v=alpha^}
\upsilon =
\hat{\alpha}(\hat{z}, p(\tau), r(\tau), r'(\tau))
\eeq
satisfies the inequality
\beq{V^'<}
    \hat{V}'(\hat{z}(\tau), p(\tau), r(\tau))
    \le- 2 \gamma
    \hat{V}(\hat{z}(\tau), p(\tau), r(\tau)).
\eeq
\item
If the functions $p, r$ are bounded and
the inequality \eq{V^'<} holds on the solutions
of the system \eq{BSCSz}, \eq{BSCSzeta}, \eq{v=alpha^},
then this system is globally asymptotically stable.
\item
If
\beq{Bb=0}
\forall p \in R^{k_p} \; \;
B(0, 0, p) = 0 \; \mbox{\rm ~and~}
b(0, 0, p) = 0,
\eeq
then
\beq{alpha^=0}
\forall p \in R^{k_p}, \; r,r_1 \in R^{k_r}
\;\; \hat{\alpha}(0, p, r, r_1) = 0.
\eeq
\item
If, in addition to listed assumptions, the function
$r'$ is bounded, then
\beq{limalpha^=0}
\lim_{\tau \to \infty}
\hat{\alpha}(\hat{z}(\tau), p(\tau), r(\tau), r'(\tau))
= 0
\eeq
holds on  the solutions of the closed-loop system
\eq{BSCSz}, \eq{BSCSzeta}, \eq{v=alpha^}.
\end{itemize}
\end{theorem}

{\bf Proof.}
The continuity of the function $\hat{\alpha}$
follows from  \eq{alpha^} and
properties of functions $\alpha, B, b, \beta,$ and $V.$
It should be noted that the function $1/\beta$
is defined and continuous
for all $z \in R^{k_z}, \; \zeta \in R, \; \tau \ge 0$
because $\beta$ is continuous and non-vanishing
and that the function $D(z, \zeta, \alpha(z, r), r)$
is defined and continuous
for all $z \in R^{k_z}, \; \zeta \in R, \; r \in R^{k_r}$
due to the continuity of the functions
$B$ and $\D{B}{\zeta}.$

Let us show that  $\hat{V}$
satisfies the conditions \eq{V=0} -- \eq{infV>0}.
The equality \eq{V=0} is obvious.

Consider the inequality \eq{V>0}.
Let $\hat{z} = \col(z, \zeta) \not= 0.$
We have
$\hat{V}(\hat{z}, p, r) > 0$ for $z \not= 0,$
since $\hat{V}(\hat{z}, p, r) \ge V(z, p, r).$
If $z = 0$ then
$\hat{V}(\hat{z}, p, r) =
\delta (\zeta - \alpha(0, r))^2/2,$
and, by virtue of \eq{alpha=0},
$\hat{V}(\hat{z}, p, r) = \zeta^2/2 > 0$
for $\zeta \not= 0.$
The inequality \eq{V>0} is proved.

We show  \eq{infV>0} by reductio ad absurdum.
Suppose  this inequality is not fulfilled.
Then there are
$\varepsilon>0,$  $\mathcal{E}>0,$
and sequences
$\{\hat{z}_\kappa\}_{\kappa = 1}^\infty,$
$\{r_\kappa\}_{\kappa = 1}^\infty,$
$\{p_\kappa\}_{\kappa = 1}^\infty,$
$(|\hat{z}_\kappa| > \varepsilon, \;
|{r}_\kappa| < \mathcal{E}, \;
|{p}_\kappa| < \mathcal{E},
\; \kappa = 1, 2, \ldots)$
such that
\beq{limV^=0}
\lim_{\kappa \to \infty}
\hat{V}(\hat{z}_{\kappa}, p_{\kappa}, r_{\kappa}) = 0.
\eeq
Let
$\hat{z}_\kappa = \col(z_\kappa, \zeta_\kappa).$
Then by virtue of \eq{limV^=0},
\beq{limV=0}
\lim_{\kappa \to \infty}
{V}({z}_{\kappa}, p_{\kappa}, r_{\kappa}) = 0.
\eeq
The limit \eq{limV=0} and the inequality \eq{infV>0} imply
\beq{limz=0}
\lim_{\kappa \to \infty} z_\kappa = 0.
\eeq
From the continuity of the function $\alpha$ and the limit \eq{limz=0}
it follows that
\beq{limalpha=0}
\lim_{\kappa \to \infty} \alpha(z_\kappa, r_\kappa) = 0.
\eeq
The substitution of \eq{limz=0} and \eq{limalpha=0}
in \eq{limV^=0} gives
$\dlim_{\kappa \to \infty} \zeta_\kappa = 0$
and, consequently,
$\dlim_{\kappa \to \infty} \hat{z}_\kappa = 0.$
The latter contradicts the assumption made about the sequence
$\{\hat{z}_\kappa\}_{\kappa = 1}^\infty.$
This contradiction
proves that the function $\hat{V}$ satisfies  \eq{infV>0}.

The  derivative of the function $\hat{V}$
along the trajectories  of the closed-loop system
\eq{BSCSz}, \eq{BSCSzeta}, \eq{v=alpha^}
has the form
\begin{eqnarray*}
    \hat{V}'(\hat{z}, p, r) =
    \D{V(z, p, r)}{z} B(z,\zeta, p) +
    \D{V(z, p, r)}{p} p' + \D{V(z, p, r)}{r} r' +\\
    \delta (\zeta - \alpha(z, r))[b(z,\zeta, p)
    + \beta(z,\zeta, p) \upsilon -
    \D{\alpha(z, r)}{z} B(z,\zeta, p) -
    \D{\alpha(z, r)}{r} r'] \\
    =
    \left\{
    \D{V(z, p, r)}{z} B(z, \alpha(z, r), p) +
    \D{V(z, p, r)}{p} p' + \D{V(z, p, r)}{r} r'
    \right\} + \\
    \left\{
    \D{V(z, p, r)}{z} B(z,\zeta, p) -
    \D{V(z, p, r)}{z} B(z, \alpha(z, r), p) +
    \right.\\ \left.
    \delta (\zeta - \alpha(z, r))[b(z,\zeta, p) +
    \beta(z,\zeta, p) \upsilon -
    \D{\alpha(z, r)}{z} B(z,\zeta, p) -
    \D{\alpha(z, r)}{r} r']
    \right\}.
\end{eqnarray*}
The expression in the first curly
braces is the derivative of $V$ along the
trajectories of the system \eq{BSCSz}, \eq{zeta=alpha}.
By virtue of \eq{alpha^},
the term in the second curly braces is equal to
$-\gamma \delta (\zeta - \alpha(z, r))^2.$
Thus, taking into account \eq{V'<},
we obtain  \eq{V^'<}.

The inequality \eq{V^'<} implies
that  the limit
\beq{V^->0}
\lim_{\tau \to \infty}
\hat{V}(\hat{z}(\tau), p(\tau), r(\tau)) = 0
\eeq
takes place on the solutions of the
system \eq{BSCSz}, \eq{BSCSzeta}, \eq{v=alpha^}.

Let us show that if
$\max(|p(\tau)|, |r(\tau)|) < \mathcal{E}$
for some $\mathcal{E}$ and all $\tau \ge 0,$
then the system
\eq{BSCSz}, \eq{BSCSzeta}, \eq{v=alpha^}
is globally asymptotically stable.
If this is not the case, then
there are $\varepsilon > 0,$  a solution $\hat{z}$
of the system \eq{BSCSz}, \eq{BSCSzeta}, \eq{v=alpha^},
 and a sequence
$\{\tau_\kappa\}_{\kappa = 1}^\infty, \;
\tau_\kappa \to_{\kappa \to \infty} \infty$
such that
\beq{infz>eps}
\inf_{\kappa = 1, 2, \ldots}
|\hat{z}(\tau_\kappa)| > \varepsilon.
\eeq
By virtue of \eq{infV>0}, we have
$\inf_{|\hat{z}| > \varepsilon, \;
\max(|p|, |r|) < \mathcal{E}}
\hat{V}(\hat{z}, p, r) > 0,$
therefore  it follows from \eq{infz>eps} that
$\hat{V}(\hat{z}(\tau_\kappa), p(\tau_\kappa), r(\tau_\kappa))
\not\to_{\kappa \to \infty} 0.$
This contradicts \eq{V^->0}.
Thus the system \eq{BSCSz}, \eq{BSCSzeta}, \eq{v=alpha^}
is globally asymptotically stable.

The equality \eq{alpha^=0} follows from \eq{alpha^}.
This implication is based on the equalities
\eq{Bb=0},
\eq{alpha=0}
and on the identity
$\D{V(0, p, r)}{z} \equiv 0,$
which follows from the fact that for any
$p$ and $r$ the function
$V(z, p, r)$ achieves minimum when $z = 0.$

The limit \eq{limalpha^=0} results
from the  continuity of the function $\hat{\alpha},$
the boundedness of the functions $r$ and $r',$ and
from the equality \eq{alpha^=0}.
\eproof \\

\noindent
{\bf 5.3.~Stabilization of $x$-subsystem.}
Consider the stabilization problem
for the $x$-subsystem\footnote{
It should be noted
that the $x$-subsystem \eq{x-CS} is not
a transformation of the kinematic model
of the Chaplygin sled,
which is third order system  and
is considered in Section~6.}
 \eq{x-CS} of the system \eq{x-CS}, \eq{s-CS}.
The feedback, that solves this problem, is used to initiate
the recursive process of designing the stabilizing control
law for the system \eq{x-CS}, \eq{s-CS}.

The inputs of the system \eq{x-CS} are $w_1$ and $\ts_1.$
Denote $E(w_1, \ts_1)$ the right-hand side
of  \eq{x-CS}.
It is evident that the equation
\beq{E=e}
E(w_1, \ts_1) = e
\eeq
is solvable for any vector $e \in R^2.$
Taking into account that
on the solutions of the closed-loop system
$\ts_1$ has to tend to zero,
we look for a solution of \eq{E=e} that satisfies the inequality
\beq{|ts_1|<}
|\ts_1| < \pi/2.
\eeq
Let us rewrite  \eq{E=e} as
\beq{=c}
w_1 \matrix{\cos \ts_1\\ \sin \ts_1} = c,
\eeq
where
$$c = G(\bs_1)^{-1} (e + \bw_1
\matrix{\cos \bs_1\\\sin \bs_1}) =
G(\bs_1)^{-1} e + \bw_1
\matrix{1\\0},$$
$$
   G(\bs_1)^{-1} =
   \matrix{\cos \bs_1 & \sin \bs_1 \\
   -\sin \bs_1 & \cos \bs_1}.
$$
For nonzero $c,$ the equation \eq{=c}
has a solution satisfying the
inequality \eq{|ts_1|<} only if $c_1 \not= 0.$
This condition can be guaranteed if the vector $e$
satisfies the inequality $|e| < 1$
since in this case ${\rm sign}~c_1 = \bw_1.$
Define
\beq{e=}
e = - \frac{\tanh(\gamma |\tx|)}{|\tx|} \tx.
\eeq
Consider the Lyapunov function candidate
\beq{V0}
    V_0(\tx)=\sinh^2(\gamma |\tx|).
\eeq
Calculating the derivative of $V_0$ along the trajectories of
the system \eq{x-CS} and assuming that right-hand side
of \eq{x-CS} is equal to $e,$
by virtue of \eq{e=} we obtain
\beq{V'=-2gV}
    V_0'(\tx)=2 \gamma \frac{\sinh(\gamma |\tx|)
    \cosh(\gamma |\tx|)}{|\tx|}
    \langle\tx,e\rangle =
    - 2 \gamma V_0(\tx).
\eeq
Thus, to guarantee the stability of the closed-loop system
it is sufficient to put
\beq{w=,s=}
\begin{array}{l}
w_1 = \bw_1 |c|
\defeq \lambda(\tx, \bs_1, \bw_1), \\
\ts_1 = \arctan(c_2/c_1)
\defeq \alpha_0(\tx, \bs_1, \bw_1),
\end{array}
\eeq
where
$$
c = -\frac{\tanh(\gamma |\tx|)}{|\tx|}
G(\bs_1)^{-1} \tx + \bw_1 \matrix{1\\0}.
$$
In such a way, we arrive at
\begin{proposition}\hspace{-1ex}{\bf :}
The closed-loop system \eq{x-CS}, \eq{w=,s=}
is globally asymptotically stable and has
the Lyapunov function \eq{V0} satisfying the inequality
\eq{V'=-2gV}.\\
\end{proposition}

\noindent
{\bf 5.4.~Recursive design of stabilizing feedback.}
Using  the feedback \eq{w=,s=},
transform the system \eq{x-CS}, \eq{s-CS}
to the form that is convenient for the recursive
application of the backstepping procedure.
Let
$
w_1= \lambda(\tx, \bs_1, \bw_1),
$
where $\lambda$ is defined by \eq{w=,s=}.
Define functions
$p^i : R \to R^{i+2}, \; i = 0, \ldots, n,$
as follows:
$p^0 = \col(\cos \bs_1, \sin \bs_1),$
$p^i(\tau) =
\col(\cos \bs_1, \sin \bs_1, \bs_2,
\ldots, \bs_{i+1}), \;
i = 1, \ldots, n-1,$
$p^{n}(\tau) =
\col(\cos \bs_1, \sin \bs_1, \bs_2,
\ldots, \bs_n, \bv_2/\bv_1).$
Then the system \eq{x-CS}, \eq{s-CS}
can be written as
\begin{eqnarray}
   && \tx' = B_0(\tx, \ts_1, p^0(\tau), \bw_1),
    \label{PCS1}\\
   && \ts'_i=b_i(\tx, p^i(\tau), \bw_1) +
    \beta_i(\tx, p^i(\tau), \bw_1) \ts_{i+1},
    \; i=1,\ldots,n,
    \label{PCS2}
\end{eqnarray}
where
\begin{eqnarray*}
   && B_0(\tx, \ts_1, p^0(\tau), \bw_1) =
    - \bw_1 \matrix{\cos \bs_1 \\ \sin \bs_1} +
    \lambda(\tx, \bs_1, \bw_1)
    G(\bs_1) \matrix{\cos \ts_1 \\ \sin \ts_1},
    \\
   && b_i(\tx, p^i(\tau), \bw_1) =
    (\lambda(\tx, \bs_1, \bw_1) - \bw_1)
    \bs_{i+1}, \
    \beta_i(\tx, p^i(\tau), \bw_1) =
    \lambda(\tx, \bs_1, \bw_1),
    \\
   && i=1,\ldots,n-1,\\
   && \bs_{n+1} = w_2, \; b_n = 0, \; \beta_n = 1.
\end{eqnarray*}

The design of the stabilizing feedback is performed by
the recursive use of the backstepping procedure to
 subsystems of the system \eq{PCS1}, \eq{PCS2},
wherein we successively increase
the number of equations in subsystems.

It is convenient to represent the $i$th
subsystem in the form
\begin{eqnarray}
 &&   (z^{i-1})'=B_{i-1}(z^{i-1}, \ts_i,
    p^{i}(\tau), \bw_1),
    \label{CSBi} \\
 &&   \ts'_i=b_i(z^{i-1}, p^{i}(\tau), \bw_1)
   + \beta_i(z^{i-1}, p^{i}(\tau), \bw_1)\ts_{i+1},
   \label{CSbi}
\end{eqnarray}
where
$z^i = \col(\tx, \ts_1, \ldots, \ts_i), \;
i=0, \ldots, n,$ is the state vector of the $i$th
subsystem,
$$
\begin{array}{l}
B_{i-1}(z^{i-1}, \ts_i, p^{i}(\tau), \bw_1) =
\col(B_0(z^{0}, \ts_1, p^{0}(\tau), \bw_1),
b_1(z^{0}, p^{1}(\tau), \bw_1)+\\
\beta_1(z^{0}, p^{1}(\tau), \bw_1)\ts_{2}, \ldots,
b_{i-1}(z^{i-2}, p^{i-1}(\tau), \bw_1) +
\beta_{i-1}(z^{i-2}, p^{i-1}(\tau), \bw_1)\ts_{i}).
\end{array}
$$
Since $\bw_1$ is constant, it is considered as a parameter.

Let us describe the $i$th step of the recursion.
Suppose that on the previous step 
functions
$\alpha_{i-1}(z^{i-1}, p^{i-1}, \bw_1)$
and
$V_{i-1}(z^{i-1}, p^{i-1}, \bw_1)$
were  constructed
such that the derivative of the function $V_{i-1}$
along the trajectories of the closed-loop system \eq{CSBi},
$$
\ts_i = \alpha_{i-1}(z^{i-1}, p^{i-1}(\tau), \bw_1)
$$
satisfies the inequality
$$
V'_{i-1}(z^{i-1}(\tau), p^{i-1}(\tau), \bw_1)
\le - 2 \gamma V_{i-1}(z^{i-1}(\tau), p^{i-1}(\tau), \bw_1).
$$
On the first step
we use the defined in Subsection~5.3
functions $\alpha_0$ and $V_0$
that satisfy the above assumption.

Choose an arbitrary $\delta_i > 0$ and define functions
$\hat{\alpha}_{i-1},$  $\hat{V}_{i-1}$
according to \eq{alpha^}, \eq{V^} with
$\alpha = \alpha_{i-1}, \; V = V_{i-1}, \;
z = z^{i-1}, \; \hat{z} = z^i, \;
r = p^{i-1}, \; r_1 = (p^{i-1})', \; p = p^i, \;
B = B_{i-1}, \; b = b_i, \; \beta = \beta_i, \; \delta = \delta_i.$
Note that using the equations \eq{x-CSdes}, \eq{s-CSdes}
we can represent the functions
$p^{i-1}$ and $(p^{i-1})'$
in terms of the function $p^i$
as follows:
\begin{eqnarray*}
& & p^{i-1} = \col(p^i_1, \ldots, p^i_{i-1})
\defeq P_i(p^i), \\
& & \begin{array}{l}
(p^{i-1})' =
\col(-\bs'_1 \sin \bs_1, \bs'_1 \cos \bs_1,
\bs'_2, \ldots, \bs'_{i+1}) = \\
 \col(- p^i_2 p^i_3 \bw_1, p^i_1 p^i_3 \bw_1,
p^i_{4} \bw_1, \ldots, p^i_i \bw_1)
\defeq P'_i(p^i, \bw_1).
\end{array}
\end{eqnarray*}
Define functions
\begin{eqnarray*}
& & \alpha_i(z^i, p^i, \bw_1) =
\hat{\alpha}_{i-1}(z^i, p^i, P_i(p^i),
P'_i(p^i, \bw_1), \bw_1), \\
& & V_i(z^i, p^i, \bw_1) =
\hat{V}_{i-1}(z^i, p^i, P_i(p^i, \bw_1), \bw_1).
\end{eqnarray*}
By virtue of Theorem~\ref{backstep}, the functions
$\alpha_i, V_i$ have the same properties as the functions
$\alpha_{i-1}, V_{i-1}.$
Consequently, the recursion can be continued.

On the $n$th step of the recursion, the function $\alpha_n$
is defined such that the system \eq{PCS1}, \eq{PCS2},
$
w_2 = \alpha_n(z^n, p^n(\tau), \bw_1)
$
is globally asymptotically stable.

Turning to the system \eq{x-CS}, \eq{s-CS},
we define the following feedback function
\beq{W=}
\Psi(\tx, \ts, \bs, \bv) =
\matrix{\lambda(\tx, \bs_1, \bw_1)\\
\alpha_n(z^n, p^n, \bw_1)},
\eeq
where
$p^n = \col(\cos \bs_1, \sin \bs_1,
\bs_2, \ldots, \bs_n, \bv_2/\bv_1), \;
z^n = \col(\tx, \ts), \;
\bw_1 = \mbox{\rm sign }\bv_1.
$
The result obtained can be formulated as
\begin{proposition}\hspace{-1ex}{\bf :}
\label{stab}
Let
the functions
$\bs_i, \; i =2, \ldots, n,$
be bounded, then the  closed-loop system
\eq{x-CS}, \eq{s-CS},
\beq{w=W}
w = \Psi(\tx, \ts, \bs, \bv),
\eeq
where $\Psi$ is defined by \eq{W=},
is globally asymptotically stable and
\beq{limw}
\lim_{\tau \to \infty} w(\tau) = \col(\bw_1, 0).
\eeq
\end{proposition}

{\bf Proof.}
To prove the proposition it is sufficient to show that the
conditions of Theorem~\ref{backstep} are fulfilled on each
step of the recursion.

Let us begin from the smoothness of the considered functions.
The functions
$\frac{\tanh(\gamma |\tx|)}{|\tx|} \tx$ and
$V_0(\tx) = \sinh^2(\gamma |\tx|)$
are infinitely differentiable for all $\tx \in R^2,$
the matrices
$G(\bs_1), G(\bs_1)^{-1}$
are infinitely differentiable for all $\bs_1 \in R.$
Therefore for fixed $\bw_1 = \pm 1$ the functions
$\lambda(\tx, \bs_1, \bw_1),$
$\alpha_0(\tx, \bs_1, \bw_1),$
$B_0(\tx, \ts_1, p^0, \bw_1),$
$b_1(\tx, p^1, \bw_1),$
$\beta_1(\tx, \bs_1, \bw_1)$
are infinitely differentiable with respect
to the other arguments.

Because of the function $\lambda$  definition
we have
$\beta_i(\tx, p^i, \bw_1) =
\lambda(\tx, \bs_1, \bw_1) \not= 0$
for all
$\tx \in R^2, \bw_1 = \pm 1, p^i \in R^{i+2},
i = 1, \ldots, n.$
Since
\beq{lambda(0)}
\lambda(0, \ts_1, \bw_1) = \bw_1
\eeq
for all
$\ts_1 \in R$ and $\bw_1 = \pm 1,$
we have from the definition of the functions
$b_i$ and $B_{i-1}$
that the equalities
$b_i(0, p^i, \bw_1) = 0$ and
$B_{i-1}(0, 0,  p^i, \bw_1) = 0$
hold for all
$i = 1, \ldots, n,$
$p^i \in R^{i+2}, \bw_1 = \pm 1.$
From \eq{w=,s=} it follows that
the equality
$
\alpha_0(0, p^0, \bw_1) = 0
$
is fulfilled for all $p^0 \in R^2.$
By virtue of \eq{V'=-2gV}, the function $V_0$
satisfies  \eq{V'<}.

In such a way,  all the conditions
of Theorem~\ref{backstep} are fulfilled on the first
and all subsequent steps during the recursion.
This implies that the closed-loop system
\eq{x-CS}, \eq{s-CS}, \eq{W=}
is globally asymptotically stable; and, moreover,
\eq{alpha^=0} implies
$\dlim_{\tau \to \infty} w_2(\tau) = 0.$
The limit
$\dlim_{\tau \to \infty} w_1(\tau) = \bw_1$
follows from \eq{lambda(0)},
the asymptotic stability of the closed-loop system
and from the uniform continuity of $\lambda.$
\eproof \\

\noindent
{\bf 5.5.~Main result.}
To obtain the control law for the system
\eq{CSx}, \eq{CSy}, the variables $\tx, \ts, w, \tau$
in the control law \eq{W=} should be transformed into
the initial variables $x, y, u, t$
and the function $p^n(\tau)$ should be expressed
in terms of  the trajectory
$q^\des(t(\tau)) = \col(\dx(t(\tau)), \dy(t(\tau)))$
and  of the input $\du(t(\tau)).$

The functions $\bx, \bs, \bv$ are defined in such a way
that for all $\tau \ge 0$ the equalities
\beq{S&F}
\begin{array}{l}
\bx(\tau) = \dx(t(\tau)), \\
\bs(\tau) = S(\dy(t(\tau))), \\
\bv(\tau) = F(\dy(t(\tau))) \du(t(\tau))
\end{array}
\eeq
hold.
The input  $v$ can be obtained
from  \eq{w}
and the last among the equalities in \eq{S&F}
\beq{w->v}
\begin{array}{l}
v(t) = |\bv_1(\tau(t))| w(\tau(t)) +
\matrix{0\\ \bv_2(\tau(t))}=
|\du_1(t)| w(\tau(t)) + F_2(\dy(t)) \du(t),
\end{array}
\eeq
where
$F_2 = \matrix{0 & 0\\0 & 1} F.$
The formulae \eq{S&F}, \eq{w->v}, and \eq{F} define
the desired feedback function
\beq{U=}
\begin{array}{l}
\Phi(q, q^\des, \du) =
~[F(y)]^{-1} \left\{|\du_1|
\Psi(x - \dx, S(y) - S(\dy), S(\dy), F(\dy) \du) +
F_2(\dy) \du \right\},
\end{array}
\eeq
where
$q, q^\des \in R^{n+2}, \;
q = \col(x, y), \;
q^\des = \col(\dx, \dy), \;
\du \in R^2$ are vectors and not functions of time.
\begin{theorem}\hspace{-1ex}{\bf :}
\label{main}Suppose that:
\begin{itemize}
\item[{\bf A.}]
The conditions of Theorem~\ref{maneuver}
are fulfilled on the manifold
$\mathcal{M} = R^2\times \mathcal{O},$
and $M$ is a maneuvering operator.
\item[{\bf B.}]
The metrics $d$ is defined on
$\mathcal{O} = R\times \check{\mathcal{O}}$
by the equation
$d(y, y') = |y_1 - y'_1| + \check{d}(\check{y}, \check{y}'),$
where
$y = \col(y_1,\check{y}), \; 
y' = \col(y'_1, \check{y}'), \; 
y_1, y'_1 \in R, \;
\check{y}, \check{y}' \in \check{\mathcal{O}},$
$\check{d}$ is a metrics on
$\check{\mathcal{O}}.$
\item[{\bf C.}]
The vectorfields $h_1,$  $h_2$
do not depend on the coordinate $y_1.$
\end{itemize}
Then the control law \eq{U=1}--\eq{u=U} with the maneuvering operator
$M$ and the feedback function $\Phi$ given by
\eq{U=} solves the  problem
of  stabilizing the distinguished point trajectories
of the system \eq{CSx}, \eq{CSy}
on the manifold $\mathcal{M}.$
In addition, on the trajectories of the closed-loop
system 
the input  $u$ is bounded.
\end{theorem}
\begin{remark}\hspace{-1ex}{\bf :}
From  \eq{U=} it follows that
for all $t \ge 0$ $\mbox{\rm sign }u_1(t) =
\mbox{\rm sign } \du_1(t)$
and, consequently, $\mbox{\rm sign }u_1(t)$
does not vary in time.
\end{remark}
\begin{remark}\hspace{-1ex}{\bf :}
The statement of Theorem~\ref{main} holds
without assumptions {\bf B} and {\bf C}
if the first component $\dy_1$
of the  trajectory $\dy$ is bounded.
However, this assumption seems to be too restrictive,
because it excludes  trajectories such as the
circle motion.
\end{remark}

{\bf Proof.}
Let
$\dx$
be a strongly admissible trajectory
of the distinguished point.
Using the maneuvering operator $M,$ define
the  trajectory $q^\des = \col(\dx, \dy)$
and the input $\du$
that correspond to $\dx$ ($(q^\des, \du) = M(\dx)$).
The change of variables from $x, y, u, t$
to $\tx, \ts, w, \tau$
transforms the system \eq{CSx}, \eq{CSy} into
the system \eq{x-CS}, \eq{s-CS}
and it transforms the
feedback 
\eq{u=Phi}
into the feedback \eq{w=W}.

By Proposition~\ref{stab},
the limits
\begin{eqnarray}
\label{tx->0}
&&\lim_{\tau \to \infty} \tx(\tau) = 0,\\
\label{ts->0}
&&\lim_{\tau \to \infty} \ts(\tau) = 0,
\end{eqnarray}
and \eq{limw} hold
on the solutions of the closed-loop system
\eq{x-CS}, \eq{s-CS}, \eq{w=W}.
The reversed change of variables in the formulae
\eq{tx->0}, \eq{ts->0}, \eq{limw}
gives  the limits \eq{x-dx->0},
\begin{eqnarray}
\label{lim(S)=0}
&&\lim_{t \to \infty} (S(y(t)) - S(\dy(t))) = 0, \\
\label{lim(F)=0}
&&\lim_{t \to \infty}
(F(y(t))u(t) - F(\dy(t))\du(t)) = 0.
\end{eqnarray}

According to conditions of Theorem~\ref{maneuver}, we have
$\mathcal{O} = R\times \check{\mathcal{O}}.$
Condition {\bf C} yields existence of maps
$\check{S} : \check{\mathcal{O}} \to R^{n-1}$
and
$\check{F} : \check{\mathcal{O}} \to GL(2),$
such that for all
$y = \col(y_1, \ldots, y_n) \in \mathcal{O},$
we have
\beq{S=checkS}
S(y) = \col(y_1, \check{S}(\check{y})),
\eeq
\beq{F=checkF}
F(y) = \check{F}(\check{y}),
\eeq
where
$\check{y} = \col(y_2, \ldots, y_n).$

Since $S$ is bijective, the map $\check{S}$
is a bijection $\check{\mathcal{O}}$ onto $R^{n-1}.$
Let us prove that $\check{S}$ is diffeomorphism.
It can be shown \cite{Isidori} that  the conditions
\eq{Liy=0}, \eq{Lky/=0} imply the relations
\beq{Jac}
|\det \mbox{\rm Jac } S(y)| =
|(L_{h_2}L_{h_1}^{n-1}y_1)|^n \not= 0.
\eeq
It follows from \eq{S=checkS} and \eq{Jac}  that
for all
${y} \in {\mathcal{O}}$
we have
$
\det \mbox{\rm Jac } \check{S}(\check{y}) =
\det \mbox{\rm Jac } S(y) \not= 0.
$
Thus
$\check{S}$ is a diffeomorphism (\cite{Shvartz}, Chapter~3, Theorem~29).

Consider the trajectory $s^\des = S(\dy).$
The formulae \eq{dy1}, \eq{ds_i=} and the fact that
the trajectory $\dx$ is strongly admissible imply
the boundedness of the functions
$s^\des_i, \; i = 2, \ldots, n.$
Denote
$\check{s}^\des = \col(s^\des_2, \ldots, s^\des_n),$
and let $\mathcal{D} \subset \check{\mathcal{O}}$
be a compact set such that
$\check{s}^\des(t) \in \mbox{\rm int }\mathcal{D}.$
Since the map $\check{S}^{-1}$ is uniformly continuous
on $\mathcal{D},$ the limit
$$
\lim_{t \to \infty}
(\check{S}(\check{y}(t)) - \check{S}(\check{y}^\des(t)))
= 0
$$
implies the limit
\beq{limcheck=0}
\lim_{t \to \infty}
\check{d}(\check{y}(t), \check{y}^\des(t)) = 0.
\eeq
The limit \eq{y-dy->0}
follows from condition {\bf B} of the theorem,
from the limits \eq{lim(S)=0}, \eq{limcheck=0}, and
the equalities $y_1 = s_1, \; y^\des_1 = s^\des_1.$

To prove  \eq{u-du->0} let us show first
the boundedness of the input function $\du$.
From the assumption that $\dx$ is strongly admissible
and from the formulae \eq{RefU1} -- \eq{dv_2=} it follows that
$v^\des$ is bounded.
By virtue of \eq{F} and \eq{F=checkF}, we have
$\du = (\check{F}(\check{S}^{-1}(\check{s}^\des)))^{-1}
v^\des.$
The latter equation implies the boundedness of $\du$
taking into account the inclusion
$\check{s}^\des(\tau) \in \mathcal{D}$ and
the boundedness of the continuous map
$(\check{F}\circ\check{S}^{-1})^{-1}$ on the compact set
$\mathcal{D}.$

Let us rewrite \eq{lim(F)=0} as follows:
\beq{limF=...}
\begin{array}{ll}
\lim_{t \to \infty} &
(F(y(t))u(t) - F(\dy(t))\du(t)) = \\
\lim_{t \to \infty} & \left\{
[\check{F}(\check{S}^{-1}(\check{s}(t)))u(t) -
\check{F}(\check{S}^{-1}(\check{s}(t)))\du(t)] +
\right. \\ & \left.
[\check{F}(\check{S}^{-1}(\check{s}(t)))\du(t) -
\check{F}(\check{S}^{-1}(\check{s}^\des(t)))\du(t)]
\right\} = 0.
\end{array}
\eeq
Considering that the map $\check{F}\circ\check{S}^{-1}$
is  continuous on $\mathcal{D}$
and that the input $\du$ is bounded,
\eq{lim(S)=0} implies
\beq{limFu-limFud=0}
\lim_{t \to \infty}
[\check{F}(\check{S}^{-1}(\check{s}(t)))\du(t) -
\check{F}(\check{S}^{-1}(\check{s}^\des(t)))\du(t)] = 0.
\eeq
The limit
\beq{limF(u-du)=0}
\lim_{t \to \infty} \{
\check{F}(\check{S}^{-1}(\check{s}(t)))[u(t) -\du(t)]
\} = 0
\eeq
follows from \eq{limF=...} and \eq{limFu-limFud=0}.
Inasmuch as $\check{s}(t) \in \mathcal{D}$
for all sufficiently large $t$
and the map $\check{F}\circ\check{S}^{-1}$ is bounded
on $\mathcal{D},$  \eq{limF(u-du)=0} implies the limit
\eq{u-du->0}.

The boundedness of the input $u$ follows from
the boundedness of $\du$ and from \eq{u-du->0}.
\eproof \\

\noindent
{\large \bf 6.~~Trajectory stabilization for a truck
with multiple trailers.}\\

\begin{figure}[t]
\begin{center}
\includegraphics[width=10cm]{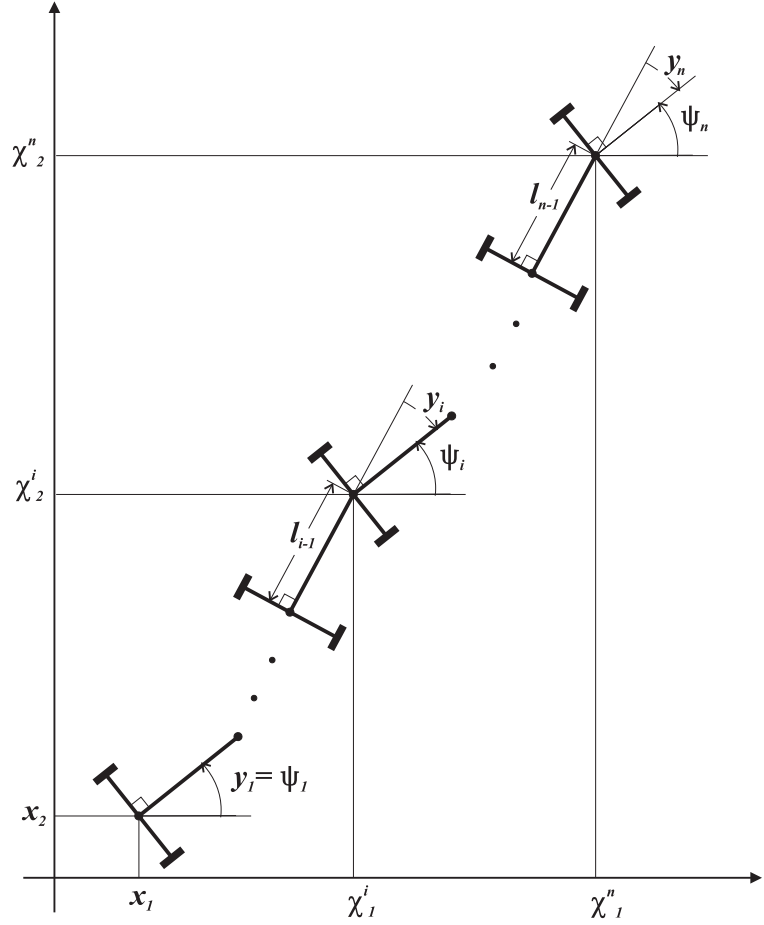}
\end{center}
\caption{Kinematic scheme of a truck with multiple
trailers.}
\label{roadtrain}
\end{figure}

Consider a wheeled system that consists of a truck
and several half-trailers; the kinematic scheme 
is shown in Fig.~\ref{roadtrain}.
Possible collisions of different parts of the vehicle
are ignored.
The configuration manifold of the system is
$\mathcal{Q} = R^3\times S^{n-1}$
and the vector of coordinates is
$q = \col(x_1, x_2, y_1, \ldots, y_n),$
where
$x_1, x_2$
are the Cartesian coordinates of the distinguished point
of the system, which is taken to be the midpoint of the axle
of the first half-trailer (trailers are enumerated
starting from the tail-end),
$y_1$ is the heading angle of this half-trailer,
$y_i$ is the angle between the axles of $i$th and
$(i-1)$th half-trailers, $i = 2 \ldots, n-2,$
$y_{n-1}$ is the angle between the axle of the last
half-trailer and the rear axle of the truck,
$y_n$ is the angle between the axles of the truck.
The kinematic model of the system has
the following form
\beq{KMT} 
\begin{array}{lll}
\dot{x}_1 & \= & u_1\cos y_1, \\
\dot{x}_2 & \= & u_1\sin y_1, \\
\dot{y}_i&\=&u_1 \eta_i(y_2,\ldots,y_{i+1}),
\; i=1,\ldots,n-1,\\
\dot{y}_n&\=&u_2,
\end{array}
\eeq 
where
$
\eta_1(y_2) \! = \! l_1^{-1} \tan y_2, \, 
\eta_i(y_2, \ldots, y_{i+1})\! = \!
(l_i^{-1}  \tan y_{i+1}
 - l_{i-1}^{-1} \sin y_i)
\prod_{k=2}^i\sec y_k, \, 
i \! = \! 2, \ldots, n\!-\!1, 
$
\newline
$u_1$ is the longitudinal velocity of
the first half-trailer,
$u_2$ is the angular velocity of the truck forward axle spin
with respect to the body of the truck.
The equations \eq{KMT} 
are derived in Appendix~B.
The system \eq{KMT} 
is a special
case of the system \eq{CSx}, \eq{CSy} with
$
h_1 = \col(\eta_1, \ldots, \eta_{n-1}, 0),\; \;
h_2 = \col(0, \ldots, 0, 1).
$

The system \eq{KMT} is defined on the manifold
$\mathcal{K} = \{ q \in \mathcal{Q} \; | \;
\cos y_i \not= 0, \; i =2, \ldots, n \}.$
The kinematic model \eq{KMT} is not
defined when any two neighbor axles are orthogonal.
When $n>1$ the manifold $\mathcal{K}$ is disconnected,
being the union
$\mathcal{K} = \bigcup_{\mu \in \mathrm{C}} \mathcal{M}_\mu$
of components $\mathcal{M}_\mu$
with the multi-index
$\mu = (\mu_1, \ldots, \mu_{n-1})$
taking on the values
among the corners $\mathcal{C}$
of the $(n-1)$-dimensional cube,
$\mathcal{C} = \{ \mu \in R^{n-1} \; | \;
\mu_i = 0 \mbox{\rm ~or~} 1, \; i =  1, \ldots, n-1\}.$
Each $\mathcal{M}_\mu, \; \mu \in \mathcal{C},$
is connected and has the form
$\mathcal{M}_\mu = \{q \in \mathcal{Q} \; | \;
\mu_{i-1} \pi - \pi/2 < y_i < \mu_{i-1} \pi + \pi/2, \;
i = 2, \ldots, n \}.$
It should be noted that for $\mu \not= 0$ the submanifold
$\mathcal{M}_\mu$ includes exotic configurations
with a neighbor half-trailers having the opposite orientation.
The manifold $\mathcal{M}_\mu$ can be represented as
$\mathcal{M}_\mu = R^2\times \mathcal{O}_\mu,$
where
$\mathcal{O}_\mu = \{y \in R\times S^{n-1} \; | \;
\mu_{i-1} \pi - \pi/2 < y_i < \mu_{i-1} \pi + \pi/2, \;
i = 2, \ldots, n \}.$
The metrics on the manifold
$\mathcal{O}_\mu, \; \mu \in \mathcal{C},$
is defined as
$d(y',y'') = \sum_{i=1}^n |y'_i-y''_i|.$
\begin{proposition}\hspace{-1ex}{\bf :}
The system \eq{KMT}
satisfies the conditions of Theorem~\ref{main}
on each manifold
$\mathcal{M}_\mu, \; \mu \in \mathcal{C}.$
\end{proposition}

The proof of Proposition~4 is given in Appendix~C.

A corollary to Proposition~4 is the maneuverability
of the system \eq{KMT}.
It means that the midpoint of the tail-end axle of the
vehicle that has $n-2$ trailers can trace
any non-stop trajectory
$x^\des \in C^{n+1}([0, \infty), R^2)$
on the plane.
This is a characteristic property of the
smooth plane curves.
It can be considered as a mechanical description
of the smoothness of a planar curve.

According to Theorem~1, on each submanifold
$\mathcal{M}_\mu, \; \mu \in \mathcal{C},$
the system \eq{KMT} has a set of
maneuvering operators.
Let us denote by $M^+_\mu$ the maneuvering operator,
which results from choosing
sign "$+$" in \eq{RefU1} and choosing the value $s^\des_1(0)$
such that the inequality $-\pi < s^\des_1(0) \le \pi$
in equation \eq{s1(0)} holds.
Similarly, by
$M^-_\mu$  we denote the maneuvering operator,
which results from
"$-$" in \eq{RefU1}
and the value $s^\des_1(0)$
satisfying $0 \le s^\des_1(0) < 2\pi.$
Then the
operator $M^+_\mu$ defines the trajectory $q^\des$
for a
forward motion of the
tail-end trailer along the desired trajectory $\dx,$
and  $M^-_\mu$ defines the trajectory $q^\des$
for a backward motion of the
tail-end trailer along the same trajectory $\dx.$

For each $\mu \in \mathcal{C}$ define
the feedback function $\Phi_\mu$
on the manifold $\mathcal{M}_\mu$
using the formula \eq{U=}.
Now for any $\mu \in \mathcal{C}$ we construct the
feedback operators $U^+_\mu$ and $U^-_\mu$
using the equations \eq{U=1}, \eq{U=2}
with $\Phi = \Phi_\mu,$
$M = M^+_\mu$ and  $M = M^-_\mu$
respectively.
According to Theorem~3,
the control laws $U^+_\mu$ and $U^-_\mu$
solve the  problem of stabilizing the
distinguished point trajectories of
the system \eq{KMT}
on the manifold $\mathcal{M}_\mu, \; \mu \in \mathcal{C}.$

Finally,
we can  design  the feedback operators
$U^+$ and $U^-$ on the manifold
$\mathcal{K} =
\bigcup_{\mu \in \mathcal{C}} \mathcal{M}_\mu,$
using the formula \eq{U=ifUi} and the operators
$U^+_\mu, U^-_\mu, \; \mu \in \mathcal{C}.$
Since the configuration manifold $\mathcal{Q}$ is the
closure of  $\mathcal{K},$
both control laws  $U^+$ and $U^-$ solve the problem
of almost global stabilization of
the distinguished point trajectories
for the considered vehicle.

It follows from Remark~2
that for all $t \ge 0$ the operator
$U^+$ defines the positive input $u_1(t) > 0,$
and $U^-$ defines the  negative input $u_1(t) < 0.$
Thus, the control law $U^+$ ensures a forward motion of
the tail-end trailer and $U^-$ ensures a backward motion of
the tail-end trailer.
The latter control law solves
intuitively harder problem of stabilizing
the road train reverse motion along the desired trajectory.

Notice that the input $u_1$
in the system \eq{KMT}
is the longitudinal velocity of the tail-end trailer.
In practice, the speed of the vehicle is controlled
by the speed of the rear-axle assembly of the truck.
Denote this alternative input $\tilde{u}_1.$
It is straightforward to show that the values $u_1$
and $\widetilde{u}_1$ satisfy the equation
$
\tilde{u}_1 = u_1 \prod_{k=2}^{n-1}\sec y_k. \;
$
Using this equation it is possible to express
the control law in terms of inputs $\tilde{u}_1$ and $u_2.$

Consider special cases of the system \eq{KMT}.

For $n = 2$  the equations \eq{KMT}
coincide with the equations of automobile \eq{3wKM}
that were considered in Section~2.
The proposed control law solves the problem of
almost global stabilization of an automobile motion
along any  non-stop
trajectory that has three bounded derivatives.

For $n = 1$ the system \eq{KMT} takes the form
\beq{WP}
\begin{array}{lll}
   \dot{x}_1&=&u_1\cos y_1, \\
   \dot{x}_2&=&u_1\sin y_1, \\
   \dot{y}_1&=&u_2.
   \end{array}
\eeq
and describes the kinematics of the Chaplygin
sled~\cite{NF72}.
Equations \eq{WP} are also used
as the kinematic model of caterpillar vehicles,
the lunar vehicle Lunohod, the experimental robot
Hilare described in~\cite{WTSML94}.
The kinematic model \eq{WP} is defined on the whole
configuration manifold of the system $\mathcal{Q} = R^3.$
Consequently, the proposed control law globally stabilizes
strongly admissible trajectories of the Chaplygin sled.\\

\noindent
{\large \bf 7.~~Simulation.}\\

For $n = 1, 2, 3, 4$ the
rectilinear motion and the circular motion
of the system \eq{KMT}
with the constructed control law 
was simulated.
The results of simulation demonstrate the efficiency
of proposed control law.
\begin{figure}[t]
\begin{center}
\includegraphics[width=12cm]{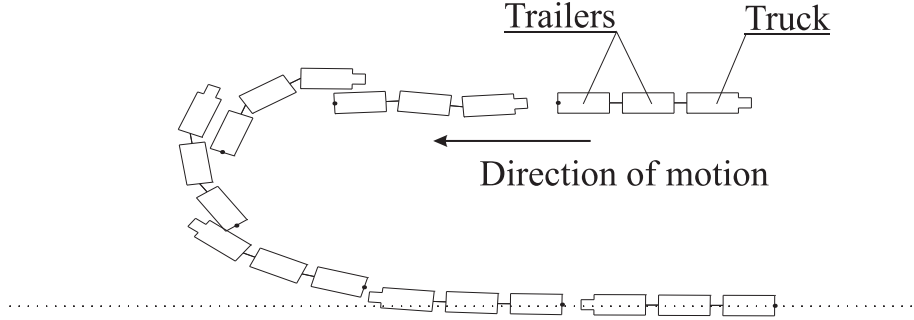}
\end{center}
\caption{U-turn of the truck pushing
two trailers in a backward direction.}
\label{Uturn}
\end{figure}
As an illustration Fig.~\ref{Uturn} shows the
sequence of vehicle positions for an U-turn of the truck pushing
two trailers in a backward direction.
The desired trajectory corresponds to the motion along the
dotted straight line.
Animated results of this and some other experiments
can be found in~\cite{GWP}.\\

\noindent
{\large \bf Appendix A. Proof of Proposition~1.}\\

By virtue of K1 the system admits
on the manifold $\mathcal{K}$
the kinematic model \eq{KM},
and due to
\newcounter{cond} \setcounter{cond}{3}
condition~\Roman{cond}
for any $q = \col(x, y) \in \mathcal{K}$
the vectorfields $g_i, i = 1,2,$
depend only on $y$-coordinates, i.e.,
$g_i(q) = g'_i(y), i = 1,2.$
Let $q = \col(x, y) \in \mathcal{K}.$
Condition K2 implies 
the equalities
$T_q(\mathcal{K}) = R^2\times T_y(\mathcal{Y})$ and
$g'_i(y) = \col(f'_i(y), h'_i(y)), \; i = 1,2,$
where
$f'_i(y) \in R^2, \; h'_i(y) \in T_y(\mathcal{Y}).$
From the nonholonomic constraint \eq{RNSC} it follows that
$f'_i(y) = \lambda_i(y) f(y), i = 1,2,$
where
$\lambda_i, i = 1,2,$ are some functions defined on
$\mathcal{Y},$ $f(y) = \col(\cos y_1, \sin y_1).$

Consider the  feedback transformation
\beq{FB1}
\matrix{\tilde{u}_1\\ \tilde{u}_2} =
\matrix{\lambda_1, & \lambda_2\\ \lambda_2, & -\lambda_1}
\matrix{u_1\\u_2}.
\eeq
Condition K3 implies that
$\lambda_1^2(y) + \lambda_2^2(y) \not= 0,$
i.e., the transformation \eq{FB1} is nonsingular for all
$y \in \mathcal{Y}.$
The transformation \eq{FB1} brings the system \eq{KM} to the
following form:
\beq{CS1}
\begin{array}{lll}
   \dot{x}  &=& \tilde{u}_1 f(y), \\
   \dot{y} &=& \tilde{u}_1 h_1(y)+ \tilde{u}_2 h_2(y),
\end{array}
\eeq
where
$h_1(y) = \lambda_1(y) h'_1(y) +
\lambda_2(y) h'_2(y), \;
h_2(y) = \lambda_2(y) h'_1(y) -
\lambda_1(y) h'_2(y).$
Equations \eq{CS1} differs from \eq{CSx}, \eq{CSy}
only by notation of inputs.\\

\noindent
{\large \bf Appendix B. Derivation of the kinematic model  for
the truck with multiple trailers.}\\

Let us introduce auxiliary variables:
$\chi_i = (\chi^1_i, \chi^2_i)$ is the vector
of the Cartesian coordinates
of the $i$th axle midpoint,
$\psi_i = \sum_{k = 1}^i y_k$ is the heading angle
of the $i$th pair of wheels,
$\tau_i = (\cos \psi_i, \sin \psi_i)$ is the unit vector,
that defines the orientation of the $i$th pair of wheels,
$\nu_i = (-\sin \psi_i, \cos \psi_i)$ is the unit vector,
that defines the orientation  of the $i$th axle,
$i = 1, \ldots, n$
(see Fig.~\ref{roadtrain}).

The nonslipping conditions for the wheels define
the nonholonomic constraints
\begin{equation}
\langle \dot{\chi}_i, \nu_i \rangle = 0,
\; \; i = 1, \ldots, n.
\label{<chi,nu>=0}
\end{equation}
In addition, the coordinates of the system
satisfy the holonomic constraints
\begin{equation}
\chi_{i+1} = \chi_i + l_i \tau_i,
\; \; i = 1, \ldots, n-1,
\label{chi=chi+ltau}
\end{equation}
that describe the articulated joints of half-trailers.
Here $l_i$ is the length of the $i$th half-trailer.

Differentiating  \eq{chi=chi+ltau}, we obtain
the equations
\begin{equation}
\dot{\chi}_{i+1} = \dot{\chi}_i +
l_i \nu_i \dot{\psi}_i,
\; \; i = 1, \ldots, n-1,
\label{chi'=chi'+lnuy'}
\end{equation}
that can be used to exclude the derivatives of the dependent
coordinates $\chi^i, i = 2, \ldots, n,$ from the
equations \eq{<chi,nu>=0}.
Thus, we deduce the equations of the nonholonomic constraints
\begin{equation}
-\sin \psi_i \dot{x}_1 + \cos \psi_i \dot{x}_2 +
\sum_{j=1}^{i-1} l_j \cos(\psi_j - \psi_i) \dot{\psi}_j = 0,
\; \; i = 1, \ldots, n,
\label{sin+cos+...=0}
\end{equation}
where $\chi^1 = \col(x_1, x_2) = x.$
The sum in the left-hand side of \eq{sin+cos+...=0} is absent
when $i = 1.$

Let us show that the nonholonomic system described by the
constraints \eq{sin+cos+...=0} admits the kinematic
model of the form \eq{CSx}, \eq{CSy}.
Scalar multiplication of  \eq{chi'=chi'+lnuy'}
by $\nu_i$ and $\tau_{i+1}$ gives the equations
\begin{equation}
\dot{\psi}_i = \upsilon_i l_i^{-1}
\tan(\psi_{i+1} - \psi_i), \;
i = 1, \ldots, n-1,
\label{y'_i=}
\end{equation}
\begin{equation}
\upsilon_{i+1} = \upsilon_i \sec(\psi_{i+1} - \psi_i), \;
i = 1, \ldots, n-1,
\label{v_i+1=}
\end{equation}
where $\upsilon_i = \langle\dot{\psi}_i, \tau_i\rangle$
is the velocity of the $i$th half-trailer.
Subject to the condition
$
\cos(\psi_{i+1} - \psi_i) \not= 0,
\; i = 1, \ldots, n-1,
$ 
the equations \eq{y'_i=}, \eq{v_i+1=}
yield the known equations of the truck with
multiple trailers kinematics~\cite{MS93}
\beq{KMTPsi}
\begin{array}{lll}
\dot{x}_1&\=&\upsilon_1\cos \psi_1, \\
\dot{x}_2&\=&\upsilon_1\sin \psi_1, \\ 
\dot{\psi}_1&\=&\upsilon_1 l_1^{-1}
\tan(\psi_{2}-\psi_1),\\
\dot{\psi}_i&\=&\upsilon_1 l_i^{-1}
\tan(\psi_{i+1}-\psi_i)
\prod_{k=2}^i\sec(\psi_k-\psi_{k-1}),
\; i=2,\ldots,n-1,\\
\dot{\psi}_n&\=&\upsilon_n,
\end{array}
\eeq
where
$\upsilon_1$ is the velocity of the tail-end half-trailer,
$\upsilon_n$ is the angular velocity of the truck
front axle spin  with respect of the truck body.
The conversion from the variables $\psi$
to the variables $y$ in
\eq{KMTPsi} gives  \eq{KMT}.\\

\noindent
{\large \bf Appendix C. Proof of Proposition~3.}\\

For fixed $\mu \in \mathcal{C}$
denote $S_i : \mathcal{O}_\mu \to R, i = 1, \ldots, n,$
the $i$th
component of the map $S,$  defined by \eq{S=}.
From the definition of the repeated
Lie derivative it follows that for $i \ge 2$
\beq{Si=sum^n}
\begin{array}{l}
S_i(y) = L_{h_1}^{i-1} y_1 =
\sum_{j = 1}^n \eta_j(y_2, \ldots, y_{j+1})
\D{}{y_j}L_{h_1}^{i-2} y_1 = 
\sum_{j = 1}^n \eta_j(y_2, \ldots, y_{j+1})
\D{S_{i-1}(y)}{y_j}.
\end{array}
\eeq
Let us show that
\beq{DjSi=0}
\D{S_i(y)}{y_j} = 0
\eeq
for  $i = 1, \ldots, n, \;
j = i+1, \ldots, n.$
To this end the equations
\beq{DSi=sum^n}
\begin{array}{l}
\D{S_i(y)}{y_j} =
\sum_{\kappa = 1}^n \D{}{y_j}
\eta_\kappa(y_2, \ldots, y_{\kappa+1})
\D{S_{i-1}(y)}{y_\kappa} + 
\sum_{\kappa = 1}^n
\eta_\kappa(y_2, \ldots, y_{\kappa+1})
\D{}{y_\kappa}\D{S_{i-1}(y)}{y_j} \;, \; i \ge 2,
\end{array}
\eeq
are used that follow from \eq{Si=sum^n}.
For $i = 1$ $S_1(y) = y_1,$
and \eq{DjSi=0} evidently holds.
Suppose that \eq{DjSi=0} is fulfilled for $i-1,$
then the equation \eq{DSi=sum^n} implies
$\D{S_i(y)}{y_j} = 0$ for $i < j \le n.$
The equalities \eq{DjSi=0} are proved.
From \eq{DjSi=0} it follows that
$S_i(y)$ does not depend on $n$ and for $i \ge 2$
\beq{Si=sum^i-1}
S_i(y) =
\sum_{j = 1}^{i-1} \eta_j(y_2, \ldots, y_{j+1})
\D{S_{i-1}(y)}{y_j}.
\eeq

The straightforward calculation gives
\beq{DiSi}
\begin{array}{l}
\D{S_i(y)}{y_{i}} =
\sum_{\kappa = 1}^{i-1} \D{}{y_{i}}
\eta_\kappa(y_2, \ldots, y_{\kappa+1})
\D{}{y_\kappa}S_{i-1}(y) + 
\sum_{\kappa = 1}^{i-1}
\eta_\kappa(y_2, \ldots, y_{\kappa+1})
\D{}{y_\kappa}\D{}{y_{i}}S_{i-1}(y) = \\
\D{}{y_{i}}
\eta_{i-1}(y_2, \ldots, y_{i})
\D{}{y_{i-1}}S_{i-1}(y) = 
l_{i-1}^{-1} \sigma_{i-1}(y_2, \ldots, y_{i-1})
\cos^{-2}(y_{i})
\D{}{y_{i-1}}S_{i-1}(y),
\end{array}
\eeq
where
$\sigma_1\equiv 1, \;
\sigma_i(y_2,\ldots,y_i)=
\prod_{k=2}^i\sec(y_k),\ i=2\ldots,n.$

From \eq{DiSi}
and the equality
$ \D{S_1(y)}{y_1} = 1$
it follows that for $i \ge 1$
\beq{DiSi=}
\D{S_i(y)}{y_{i}} =
\sigma_{i}(y_2, \ldots, y_{i})^2
\prod_{j=1}^{i-1} \left[
l_j^{-1} \sigma_j(y_2, \ldots, y_j)\right].
\eeq
The equations \eq{Si=sum^i-1} and \eq{DiSi=} give for $i \ge 2$
\beq{Si=}
S_i(y) =
\theta_i(y_2, \ldots, y_{i-1}) \tan(y_i) +
\xi_i(y_2, \ldots, y_{i-1}),
\eeq
where
$$
\begin{array}{l}
\theta_i(y_2, \ldots, y_{i-1}) =
\sigma_{i-1}(y_2, \ldots, y_{i-1})^2
\prod_{j=1}^{i-1} \left[
l_j^{-1} \sigma_j(y_2, \ldots, y_j) \right],\\
\xi_i(y_2, \ldots, y_{i-1}) =
- l_{i-1}^{-1} \sigma_{i-1}(y_2, \ldots, y_{i-1})
\tan(y_{i-1}) + 
\sum_{j = 1}^{i-2} \eta_j(y_2, \ldots, y_{j+1})
\D{S_{i-1}(y)}{y_j}.
\end{array}
$$
From the definition of $S_i$ we have
$L_{h_2}L_{h_1}^i y_1 = \D{S_{i+1}}{y_n}.$
This equality and the equations \eq{DjSi=0} and \eq{DiSi=}
imply that  \eq{Liy=0} and \eq{Lky/=0}
are fulfilled for all $y \in \mathcal{O}_\mu.$

Let us show that $S$ maps bijectively
$\mathcal{O}_\mu$ onto $R^n.$
Choose arbitrary
$s^* \in R^n$
and consider the equation
\beq{s=S(y)}
s^* = S(y).
\eeq
Define the vector $y^* \in \mathcal{O}_\mu$
by the recursive formulae
$$
\begin{array}{l}
y^*_1 = s^*_1, \;\\
y^*_i = \arctan(
\frac{s^*_i - \xi_i(y^*_1, \ldots, y^*_{i-1})}
{\theta_i(y^*_1, \ldots, y^*_{i-1})}) + \mu_{i-1} \pi, \;
i = 2, \ldots, n.
\end{array}
$$
From \eq{Si=} it is evidently follows that $y^*$
is a unique solution of \eq{s=S(y)}.

We proved that the system \eq{KMT}
satisfies the conditions of Theorem~\ref{maneuver}.

Consider condition {\bf B}  of Theorem~\ref{main}.
Manifold $\mathcal{O}_\mu$ can be represented as
$\mathcal{O}_\mu =
R\times \check{\mathcal{O}}_\mu,$
where
$\check{\mathcal{O}}_\mu = \{\check{y} \in S^{n-1} \;
| \; \mu_i \pi - \pi/2 < \check{y}_i < \mu_i \pi + \pi/2, \;
i = 1, \ldots, n-1 \}.$
On the manifold
$\check{\mathcal{O}}_\mu$
define the metrics
$\check{d}(\check{y}',\check{y}'') =
\sum_{i=1}^{n-1} |\check{y}'_i-\check{y}''_i|,$
then the metrics $d$ satisfies
condition {\bf B} of the theorem.

Condition  {\bf C} of the theorem evidently satisfied
by the definition of
the vectorfields $h_1$ and $h_2.$\\

\noindent
{\large \bf Acknowledgment}\\

The authors would like to thank B.~D.~Lubachevsky
and I.~V.~Burkov for their thoughtful comments.
\def\TrAC{{\sl IEEE Trans. Automat. Contr.}}
\def\SCL{{\sl Syst. Control Lett.}}
\def\MCSS{{\sl Mathematics of Control, Signals, and Systems}}
\def\ecc2{in {\sl Proc. 2nd European Control Conference}}
\def\cdc{in {\sl Proc. Conf. Decision Control}}
\def\cra{in {\sl Proc. IEEE Intern. Conference
on Robotics and Automation}}

\end{document}